\documentclass[11pt,twoside]{amsart}
\usepackage{geometry}
\allowdisplaybreaks

\usepackage{color}
\usepackage{amssymb}
\usepackage{dsfont}
\usepackage{amssymb,amsmath,cancel}
\usepackage{amsfonts}
\usepackage{amsthm,amssymb}
\usepackage{latexsym,amsmath,amscd}
\usepackage{enumerate}
\usepackage{color}

\numberwithin{equation}{section}

\def\proof{\noindent {\it Proof. $\, $}}
\def\endproof{\hfill $\Box$ \vskip 5 pt}

\newcommand{\be}{\begin{equation}}
\newcommand{\ee}{\end{equation}}
\newcommand{\bde}{\begin{displaymath}}
\newcommand{\ede}{\end{displaymath}}

\newenvironment{Remark}{\brem}{\erem}

\newtheorem{theorem}{Theorem}[section]
\newtheorem{lem}{Lemma}[section]
\newtheorem{pro}{Proposition}[section]
\newtheorem{cor}{Corollary}[section]
\newtheorem{rem}{Remark}[section]
\newtheorem{rems}{Remarks}[section]

\newtheorem{ex}{Example}[section]
\newtheorem{defi}{Definition}[section]
\newtheorem{hyp}{Assumption}[section]
\newtheorem{con}{Conjecture}[section]

\newcommand{\bt}{\begin{theorem}}
\newcommand{\et}{\end{theorem}}
\newcommand{\bl}{\begin{lem}}
\newcommand{\el}{\end{lem}}
\newcommand{\bp}{\begin{pro}}
\newcommand{\ep}{\end{pro}}
\newcommand{\bcor}{\begin{cor}}
\newcommand{\ecor}{\end{cor}}
\newcommand{\bcon }{\begin{con} \rm }
\newcommand{\econ }{\end{con}}
\newcommand{\lab }{\label }
\newcommand{\bd}{\begin{defi} \rm }
\newcommand{\ed}{\end{defi}}
\newcommand{\brem }{\begin{rem} \rm }
\newcommand{\erem }{\end{rem}}
\newcommand{\brems }{\begin{rems} \rm }
\newcommand{\erems }{\end{rems}}
\newcommand{\bhyp }{\begin{hyp} \rm }
\newcommand{\ehyp }{\end{hyp}}
\newcommand{\bex}{\begin{ex} \rm }
\newcommand{\eex}{\end{ex}}

\def\vert{\,|\,}
\def\Vert{\,\Big|\,}
\newcommand{\ind}{\mathds{1}}

\def\myh{A}
\def\Classmp{{\mathcal M}^p(\gg ,\tau )}
\def\Classm{{\mathcal M}^o (\gg ,\tau )}
\def\Classh{{\mathcal H}^o (\ff)}
\def\Classhp{{\mathcal H}^p (\ff)}

\def\Bpff{B^{p,\ff}}
\def\wtBpff{\widetilde{B}^{p,\ff}}

\def\F{{\mathcal F}}
\def\G{{\mathcal G}}
\def\H{{\mathcal H}}

\def\ff{{\mathbb F}}
\def\gg{{\mathbb G}}
\def\hh{{\mathbb H}}
\def\rr{{\mathbb R}}

\def\P{{\mathbb P}}
\def\E{\mathbb{E}}

\def\cro#1{\langle #1\rangle}

\def\wh{\widehat }
\def\wt{\widetilde }
\def\mz{\mu}                   %% partie $\ff$-martingale de $Z$
\def\mx{\mu^h}                 %% partie $-ff$-martingale de $X$
\def\mt{M^{(\tau)}}            %% compense de H
\def\ap{A^p}                   %% partie previsible croissante de Z
\def\ao{A^o}
\def\bmz{{\overline \mz}}
\def\bmx{{\overline{\mu}^h }}

\def\myh{A}

\def\myu{U}

\newcommand{\myeta}{\eta}
\newcommand{\myzeta}{\zeta}

\def\lab{\label}

\newcommand{\Keywords}[1]{\par\noindent{\small{\bf Keywords\/}: #1}}
\newcommand{\Class}[1]{\par\noindent{\small{\bf Mathematics Subjects Classification (2010)\/}: #1}}

\title[Integral representations for progressive enlargements]{INTEGRAL REPRESENTATIONS OF MARTINGALES FOR PROGRESSIVE ENLARGEMENTS OF FILTRATIONS}

\author[A. Aksamit, M. Jeanblanc and M. Rutkowski]{Anna Aksamit, Monique Jeanblanc and Marek Rutkowski\\
}

\address{Anna Aksamit\\
School of Mathematics and Statistics \\ University of Sydney
\\ Sydney, NSW 2006, Australia}
\email[Corresponding author]{anna.aksamit@sydney.edu.au}
\address{Monique Jeanblanc\\
Laboratoire de Math\'ematiques et Mod\'elisation d'\'Evry (LaMME), UMR CNRS 8071\\
Universit\'e d'\'Evry-Val-d'Essonne\\
Universit\'e Paris Saclay\\
23 Boulevard de France\\
91037 \'Evry cedex\\
France}
\email{monique.jeanblanc@univ-evry.fr}

\address{Marek Rutkowski
\\ School of Mathematics and Statistics \\ University of Sydney
\\ Sydney, NSW 2006, Australia  and   Faculty of Mathematics and Information Science
\\ Warsaw University of Technology
\\ 00-661 Warszawa, Poland}
\email{marek.rutkowski@sydney.edu.au}

\date{May 21, 2018}

\begin{document}
\maketitle
\vskip 20 pt

\begin{abstract}
We work in the setting of  the progressive enlargement $\gg$ of a reference filtration $\ff$ through the observation of a random time $\tau$. We study  an integral representation property for some classes of $\gg$-martingales stopped at $\tau$. In the first part, we focus on the case where $\ff$ is a Poisson filtration and we establish a predictable representation property with respect to three $\gg$-martingales. In the second part, we relax the assumption that $\ff$ is a Poisson filtration and we assume that $\tau$ is an $\ff$-pseudo-stopping time. We establish integral representations with respect to some $\gg$-martingales built from $\ff$-martingales and, under additional hypotheses, we obtain a predictable representation property with respect to two $\gg$-martingales.
\vskip 10 pt
\vskip 10 pt
\Class{60H99}
\end{abstract}

%\tableofcontents
%\newpage

%%%%%%%%%%%%%%%%%%%%%%%%%%%%%%%%%%%%%%%%%%%%%%%%%%%%%%%%%%%%%%%%%%%%%%%%%%%
%%%%%%%%%%%%%%%%%%%%%%%%%%%%%%%%%%%%%%%%%%%%%%%%%%%%%%%%%%%%%%%%%%%%%%%%%%%
\section{Introduction}  \lab{sec1}
%%%%%%%%%%%%%%%%%%%%%%%%%%%%%%%%%%%%%%%%%%%%%%%%%%%%%%%%%%%%%%%%%%%%%%%%%%%
%%%%%%%%%%%%%%%%%%%%%%%%%%%%%%%%%%%%%%%%%%%%%%%%%%%%%%%%%%%%%%%%%%%%%%%%%%%

We are interested in the stability of the predictable  representation
property in a filtration enlargement setting: under the postulate that the predictable
representation property holds for $\ff$ and $\gg$ is an
enlargement of  $\ff$, the goal is to find under which conditions the predictable
representation property holds for $\gg$.  We focus on the
progressive enlargement $\gg = ({\mathcal G}_t)_{t \in \rr_+}$ of
a reference filtration $\ff= ({\mathcal F}_t)_{t \in \rr_+}$
through the observation of the occurrence of a random time $\tau$
and we assume that the \emph{hypothesis} ($H'$) is satisfied, i.e., any
$\ff$-martingale is a $\gg$-semimartingale.
It is worth noting that no general result on the existence of a
$\gg$-semimartingale decomposition after $\tau $ of an
$\ff$-martingale is available in the existing literature, while,
for any $\tau$,  any $\ff$-martingale stopped at time $\tau$ is a
$\gg$-semimartingale with an explicit decomposition depending on
the Az\'ema supermartingale of $\tau$. Moreover, most papers in
this area are devoted to either the case of {\it honest times}
(see, e.g., Barlow  \cite{ba:sf} and  Jeulin \& Yor \cite{JeuYor,JY2}) or the case where
the {\it Jacod's absolute continuity hypothesis} holds (see, e.g., Jeanblanc \& Le
Cam \cite{jc}). The particular case where the  {\it Jacod's equivalence hypothesis} holds  is presented in  Callegaro et al. \cite{CJZ}.  More information can be found in \cite{aj}.
The stability of predictable representation  property  was studied
first by Kusuoka \cite{KS}, under the assumptions that the
filtration $\ff$ is generated by a Brownian motion,  the intensity rate of $\tau$ exists and the
filtration $\ff$ is {\it immersed} in its progressive enlargement
$\gg$, that is, any $\ff$-martingale is a $\gg$-martingale.  The \emph{immersion property}  is also known as the {\it hypothesis} ($H$) (see Br\'emaud \& Yor \cite{by}). Kusuoka  has shown
 that any $\gg$-martingale can be represented
as a sum of a stochastic integral with respect to  the Brownian
motion and a stochastic integral with respect to the compensated
$\gg$-martingale $\mt$ of the \emph{indicator process} $A:= \ind_{[\![\tau,\infty[\![}$. Under the hypotheses that the filtration $\ff$ enjoys the predictable  representation
property, $\ff$ is immersed in $\gg$, and $\tau$ avoids stopping times, this result was extended for a positive Az\'ema supermartingale $Z$   in Aksamit \& Jeanblanc \cite[Th. 3.13]{aj} and, for the case of an arbitrary supermartingale $Z$, in  Coculescu et al.  \cite{cjn}.  Barlow \cite{ba:sf} started the study  of the predictable  representation
property in the case where $\tau$ is honest. He solved the problem  for honest times avoiding $\ff$-stopping times when all $\ff$-martingales are continuous in an unpublished manuscript \cite{ba:mr}.

The predictable representation property  in a progressively
enlarged filtration  was then  studied in Jeanblanc \& Song
\cite{js3} in full generality, assuming that the predictable
representation property holds in the filtration $\ff$, with
respect to a (possibly multidimensional) martingale $M$ and that the hypothesis ($H'$) between $\ff$ and $\gg$ holds. Theorem 3.2 in \cite{js3} gives conditions for the validity of the predictable
representation property  for $\gg$-martingales with respect to the
pair $(\widehat M,\mt)$  where $\widehat M $ is the $\gg$-martingale
part of the $\gg$-semimartingale $M$. In particular, they  work under the assumption that $\G_\tau= \G_{\tau-}$ (see \eqref{htauminus} for the definition of these $\sigma$-fields). Comparing to their result we
limit ourselves to some particular set-ups but we attempt here to
derive explicit expressions for the predictable integrands, rather
than merely to establish their existence and uniqueness and we do not make the assumption that $\G_\tau= \G_{\tau-}$. Our main result in this spirit, concerning a Poisson filtration, is stated in Theorem \ref{pox2} and its further simplifications are given in Theorem \ref{cor:posi} and Corollaries \ref{cor:one} and \ref{corss}. We discuss the condition $\G_\tau=\G_{\tau-}$ in Remark \ref{rem:taumi}. Moreover, we derive different integral representations with not necessarily
predictable integrands in the specific case where $\tau$ is assumed to be a
pseudo-stopping time. The results of this type are collected in Propositions \ref{pro1n}, \ref{prox1} and \ref{pro1}.

The paper is structured as follows. In Section \ref{sec2}, we introduce the setup and notation. We also review very
briefly the classic results regarding the $\gg$-semimartingale decomposition of $\ff$-martingales stopped at $\tau $.

In Section \ref{sec3}, we assume that the filtration $\ff$ is generated by a Poisson process  $N$ with compensated
$\ff$-martingale $M$ and we consider the filtration $\gg$, which is the progressive enlargement  of $\ff$ with an arbitrary random time. Since we work with a Poisson filtration, all $\ff$-martingales are necessarily processes of finite variation and thus any $\ff$-martingale is manifestly a $\gg$-semimartingale; in other words, the \emph{hypothesis} ($H'$) is satisfied.
In Section \ref{sec3.1}, we establish the predictable representation property for a class of $\gg$-martingales stopped at time $\tau$ in terms of three martingales:  the compensated $\gg$-martingale  $\mt$ of the indicator process $A$, the
$\gg$-martingale part  $\widehat M$ of the $\gg$-semimartingale
$M$ and the $\gg$-martingale, denoted by $M^*$, which takes care of the common jumps of $\mt$ and
$\widehat M$. In Section \ref{sec3.2}, we  illustrate the above result with  a
particular example obtained with  a Cox construction.

In Section \ref{sec4}, we consider a general filtration $\ff$ and
we work under the postulate that $\tau $ is an $\ff$-{\it
pseudo-stopping time} (see Nikeghbali \& Yor \cite{NY}), so that
its Az\'ema supermartingale is an $\ff$-optional, decreasing
process. We first study in Section \ref{sec4.1} the particular case of a
random time $\tau $ independent of the filtration $\ff$. In Section
\ref{sec4.2} we consider some classes of   $\gg$-martingales
stopped at $\tau$. We derive several alternative integral
representations with respect to martingales built from
$\ff$-martingales,  under the assumption that the Az\'ema
supermartingale of $\tau $ is strictly positive. In Section
\ref{sec4.3}, the positivity assumption is  relaxed and we obtain
a more general representation result. Under additional hypotheses
on the  common jumps of specific $\ff$-martingales and $\mt$, we
obtain a predictable representation property with respect to the
two martingales $\widehat M$ and $\mt$.

\vskip 5 pt
\noindent{\bf Conventions.} We use throughout the notational convention that $\int_{t}^s = \int_{(t,s]}$. Moreover, all integrals used in what follows are implicitly postulated to be well defined without further explicit mentioning. In several (but not all) instances the existence of integrals is easy to check, since any process $h$ from the class  $\Classm$ (or $\Classmp$) is bounded
(see Section \ref{sec2.3}). Finally, we write $X=0$ when the process $X$ is indistinguishable from the null process.

%\newpage

%%%%%%%%%%%%%%%%%%%%%%%%%%%%%%%%%%%%%%%%%%%%%%%%%%%%%%%%%%%%%%%%%%%%%%%%%%%%%%%%%
%%%%%%%%%%%%%%%%%%%%%%%%%%%%%%%%%%%%%%%%%%%%%%%%%%%%%%%%%%%%%%%%%%%%%%%%%%%%%%%%%
\section{Preliminaries} \lab{sec2}
%%%%%%%%%%%%%%%%%%%%%%%%%%%%%%%%%%%%%%%%%%%%%%%%%%%%%%%%%%%%%%%%%%%%%%%%%%%%%%%%%
%%%%%%%%%%%%%%%%%%%%%%%%%%%%%%%%%%%%%%%%%%%%%%%%%%%%%%%%%%%%%%%%%%%%%%%%%%%%%%%%%

We first introduce the notation for an abstract setup in which a
reference filtration $\ff$ is progressively enlarged through
observation of a random time. Let $(\Omega, \G, \ff,\P)$ be a
probability space where $\ff$ is an arbitrary filtration
satisfying the usual conditions of right-continuity and $\P$-completeness, and such that $\F_0$ is a trivial $\sigma$-field and $\F_\infty \subseteq
\G$. Let $\tau$ be a {\it random time}, that is, a $[0,\infty]$-valued random variable on $(\Omega, \G,  \P)$.
The {\it indicator process} $A := \ind_{[\![ \tau , \infty [\![}$ is an
increasing process, which is $\ff$-adapted if and only if $\tau$ is an $\ff$-stopping time.
By convention, we shall put $A_{0-}:=0$ whenever we consider the jump of $A$ at 0.

We denote by $\gg$ the {\it progressive
enlargement} of the filtration $\ff$ with the random time $\tau$,
that is, the smallest right-continuous and $\P$-completed
filtration $\gg$ such that the inclusion $\ff \subset \gg$ holds and
the process $A$ is $\gg$-adapted so that the random time $\tau $ is a
$\gg$-stopping time. In other words, the progressive enlargement
is the smallest extension of the filtration $\ff$, which satisfies the usual conditions and
renders $\tau $ a stopping time.

%%%%%%%%%%%%%%%%%%%%%%%%%%%%%%%%%%%%%%%%%%%%%%%%%%%%%%%%%%%%%%%%%%%%%%%%%%%%%%%%
\subsection{Az\'ema Supermartingale} \label{sec2.1}
%%%%%%%%%%%%%%%%%%%%%%%%%%%%%%%%%%%%%%%%%%%%%%%%%%%%%%%%%%%%%%%%%%%%%%%%%%%%%%%%

For a filtration $\hh \in\{ \ff, \gg\}$,  we denote by $B^{p,\hh }$ (resp., $B^{o,\hh }$) the dual $\hh$-predictable (resp., the dual $\hh$-optional) projection of a process $B$ of finite variation, whereas ${}^{p,\hh} U$  (resp., ${}^{o,\hh} U$) stands for the $\hh$-predictable (resp., $\hh$-optional) projection of a process $U$. The $\hh$-predictable covariation
process of two $\hh$-semimartingales $X$ and $Y$ is denoted by $\cro{X,Y}^\hh$. Recall that
$\cro{X,Y}^\hh$ is defined as the dual $\hh$-predictable projection of the covariation process $[X,Y]$,
that is, $\cro{X,Y}^\hh := [X,Y]^{p,\hh}$.  The following definition is due to Az\'ema \cite{az}.

\bd  The c\`adl\`ag, bounded $\ff$-supermartingale $Z$ given by
$$
Z_t :=  \P(\tau>t \vert \F_t) = {}^{o,\ff} \big( \ind_{[\![ 0, \tau [\![} \big)_t = 1- {\,}^{o,\ff\!}A _t
$$\,
is called the {\it Az\'ema supermartingale} of the random time $\tau $.
\ed

Note  that $Z_{\infty }:= \lim_{t \to \infty } Z_t =\P(\tau=\infty\,|\,\F_\infty)$. By convention, we set $Z_{0-}:=1$.
It is known that the processes $Z$ and $Z_{-}$ do not vanish before $\tau $. Specifically, the random set $\{Z_{-}=0\}$ is disjoint from $[\![0,\tau]\!]$ so that the set $\{ Z = 0 \}$ is disjoint from $[\![0,\tau [\![$ (see, e.g., \cite{Y2}).

The  Az\'ema supermartingale $Z$ is of class (D) and thus it admits a (unique) {\it Doob-Meyer
decomposition} $Z=\mz -\ap $ where  $\mz$ is a positive BMO $\ff$-martingale
and $\ap$ is an $\ff$-predictable, increasing process.    Specifically, $\ap $  is the dual $\ff$-predictable projection of
the increasing process $A$, that is, $\ap = A^{p,\ff}$, and the positive
$\ff$-martingale $\mz$ is given by the equality $\mz_t:=
\E(\ap_\infty  +\ind_{\{\tau=\infty\}}\vert \F_t)$. Then $\ap_0= \E( A_0 \vert \F_0)=\P(\tau=0 \vert \F_0)=  1-Z_0$ as in \cite[Proposition 1.24]{aj} and we obtain $\mu_0=1$.  It is well known
\cite[p. 64]{Jeu}) that
the process $\mt$ given by
\begin{equation} \label{mtau}
\mt_t := A_t-\int_0^{t\wedge \tau} \frac{d\ap_s}{Z_{s-}}
\end{equation}
is a uniformly integrable $\gg$-martingale. We shall also use the dual $\ff$-optional projection $\ao := A^{o,\ff}$ of
$A$ and the positive BMO $\ff$-martingale $m$ given by $m_t := \E (\ao_{\infty }+\ind_{\{\tau=\infty\}} \vert \F_t)$. Then $  \ao_0= \E( A_0 \vert \F_0)$ and $m_0 =1$. It is also well known that $Z = m - \ao$ (see, e.g.,  \cite[Proposition 1.46]{aj}).

\brem \label{rempb}
 It is well known (see, e.g., Lemma 8.13 in Nikeghbali \cite{nk} or \cite[Lemma 1.48]{aj}) that if either: \\
(a)  $\tau$ avoids all $\ff$-stopping times, that is, $\P (\tau = \sigma<\infty )=0$ for any
$\ff$-stopping time $\sigma $ or \\
(b) all $\ff$-martingales are continuous so that all $\ff$-optional processes are $\ff$-predictable, \\
then the equalities
$\ao = \ap $ and $\mz = m$ are valid.
These properties motivated other authors to focus on models in which either the avoidance or the continuity properties are satisfied.
\erem

%%%%%%%%%%%%%%%%%%%%%%%%%%%%%%%%%%%%%%%%%%%%%%%%%%%%%%%%%%%%%%%%%%%%%%%%%%%%%%%%
\subsection{Semimartingale Decompositions for Progressive Enlargements} \label{sec2.2}
%%%%%%%%%%%%%%%%%%%%%%%%%%%%%%%%%%%%%%%%%%%%%%%%%%%%%%%%%%%%%%%%%%%%%%%%%%%%%%%%

For the reader's convenience, we recall some useful results for the
progressive enlargement of a filtration $\ff$ and an
arbitrary random time $\tau $. As was already mentioned,  no general
result giving tractable  conditions   for the stability of space of semimartingales under enlargement of filtration is known  (one can see Jeulin \cite[th. 2.6]{Jeu}) nor furnishing a $\gg$-semimartingale decomposition after $\tau
$ of an $\ff$-martingale is available in the  existing literature,
although some partial results were established for particular classes of
random times, such as: the {\it honest times} (see Barlow \cite{ba:sf}  and Jeulin \& Yor
\cite{JY2}), the random times satisfying the {\it density
hypothesis} (see Jeanblanc \& Le Cam \cite{jc}), as well as for
some random times obtained through various extensions
of the so-called {\it Cox construction} of a random time (see,
e.g., \cite{GZ,js1,js2,lr1,lr2,NY}) and recently for a class of thin random times (see Aksamit et al. \cite{acj_thin}).  In this work, we will use semimartingale decompositions
for $\ff$-martingales stopped at $\tau $ and thus we first quote some classic results regarding this
case. For the result stated in Proposition \ref{proq}, the reader is
referred to Jeulin  \cite[Proposition 4.16]{Jeu} and Jeulin \& Yor \cite[Th\'eor\`eme 1, pp.
87--88]{JeuYor}).

\bp \label{proq}
For any $\ff$-local martingale $X$, the process
\be \label{j1}
\widehat X_t := X_{t \wedge \tau}-\int _0^{t\wedge \tau} \frac {1}{Z_{s-}}\, d \cro {X,m}^\ff_s=X_{t\wedge \tau} -\int_0^{t\wedge \tau}\ind_{\{Z_{s-}<1\}} \frac
{1}{Z_{s-}} \left ( d\cro{X,\mz }^\ff_s + d  J_s\right)
\ee
where $J := \big( A\Delta X_\tau  \big)^{p,\ff }\!,$ is a $\gg$-local martingale (stopped at $\tau$).
\ep

By comparing the two decompositions in \eqref{j1}, we obtain the following well known result, which will
be used in the proof of the main result of this paper (see Theorem \ref{pox2}).

 \bcor \label{cor1}
The following equality holds for any $\ff$-local martingale $X$
\begin{equation} \label{equa1}
\int_0^{t\wedge \tau}\frac {1}{Z_{s-}}\, d  J_s = \int_0^{t\wedge
\tau}\frac {1}{Z_{s-}}\, d\cro{X,m-\mz }^\ff_s  .
\end{equation}
\ecor

\proof
For completeness, we give the proof of the corollary. We note that the processes
$$
\int_0^{\cdot\wedge \tau}\frac {1}{Z_{s-}}\, d \cro {X,m}^\ff_s
\quad \quad \textrm{and} \quad \quad
\int_0^{\cdot\wedge \tau} \ind_{\{Z_{s-}<1\}} \frac {1}{Z_{s-}}
\left ( d\cro{X,\mz }^\ff_s + dJ_s\right)
$$
are $\gg$-predictable. Hence equations \eqref{j1}  yield two equivalent forms of the Doob-Meyer
decomposition of the special $\gg$-semimartingale
$X_{\cdot \wedge \tau}$. Due to the uniqueness of the Doob-Meyer decomposition, we obtain
\begin{align*}
0&=\int_0^{t\wedge \tau}\ind_{\{Z_{s-}<1\}}  \frac
{1}{Z_{s-}}\left ( d\cro{X,\mz  }^\ff_s + dJ_s\right)
-\int _0^{t\wedge \tau}\frac {1}{Z_{s-}}\, d \cro {X,m}^\ff_s \\
&=\int_0^{t\wedge \tau}\frac {1}{Z_{s-}}\left ( d\cro{X,\mz -m
}^\ff_s + dJ_s\right)
-\int_0^{t\wedge \tau}\ind_{\{Z_{s-}=1\}} \left ( d\cro{X,\mz  }^\ff_s + dJ_s\right)\\
&=\int_0^{t\wedge \tau}\frac {1}{Z_{s-}}\left ( d\cro{X,\mz -m
}^\ff_s + dJ_s\right)
\end{align*}
where the last equality follows from Lemme 4(b) in \cite{JeuYor}. Hence equality \eqref{equa1} is valid.
\endproof

%%%%%%%%%%%%%%%%%%%%%%%%%%%%%%%%%%%%%%%%%%%%%%%%%%%%%%%%%%%%%%%%%%%%%%%%%%%%%%%%%
\subsection{Problem Formulation} \lab{sec2.3}
%%%%%%%%%%%%%%%%%%%%%%%%%%%%%%%%%%%%%%%%%%%%%%%%%%%%%%%%%%%%%%%%%%%%%%%%%%%%%%%%%

For a filtration $\hh\in\{ \ff, \gg\}$ and a random time $\tau$, we introduce two $\sigma$-fields $\H_\tau$ and $\H_{\tau-}$ defined as  \begin{equation} \label{htauminus} \left\{ \begin{array}{lcr}
 \H_\tau= \sigma \big\{Y_\tau \ind_{\{\tau<\infty\}},\, Y  \mbox{~is~an~}\,{\hh{\mbox-{optional~}}} \mbox{process} \big\}\\
\H_{\tau-}= \sigma \big\{Y_\tau \ind_{\{\tau<\infty\}},\, Y
\mbox{~is~an~}\,{\hh{\mbox-{predictable~}}} \mbox{process}\big\}.\end{array}\right. \end{equation}
Obviously, $\H_{\tau-}\subset \H_{\tau}$ and $\tau\in \H_{\tau-}$.   We recall that,  in a progressive enlargement setting, for  any $\gg$-predictable process $Y$, there exists an $\ff$-predictable process $y$ such that $Y_t\ind_{\{t\leq \tau \}}=y_t\ind_{\{t\leq \tau \}}$ (see, e.g., \cite[Proposition 2.11]{aj}) and thus we have that $\G_{\tau-}=\F_{\tau-}$. In contrast, the $\sigma$-fields $\G_\tau$ and $\F_\tau$ may differ, in general (for counterexamples, see Barlow  \cite[Page 319]{ba:sf} or \cite[Example 5.13]{aj}). Note that since the random time $\tau$ is a $\gg$-stopping time, the definitions of $\G_\tau$ and $\G_{\tau-}$ given above coincide with the usual ones,
specifically,
$$
\G_\tau =\{G \in \G_\infty \,: \, G \cap \{\tau \leq t\}\in \G_t, \,\forall t\}
$$
whereas $\G_{\tau -}$ is the smallest $\sigma$-field which contains $\G_0$ and all the sets of the form $G\cap \{t<\tau \}$ where $t>0$ and $G\in \G_t$.

Let us introduce the following notation for the classes of $\gg$-martingales studied in this work
$$
\Classm  := \big\{ Y^h_t :=  \E(h_\tau \vert \G_t)\! : h \in \Classh \big\}
$$
where $\Classh := \{ h : h \text{ a bounded $\ff$-optional process with limit at} +\infty\}$.   We denote by
$\Classhp$ the set of all processes $h$ in $\Classh$ that are $\ff$-predictable and we set
$$
\Classmp  := \big\{ Y^h_t :=  \E(h_\tau \vert \G_t)\! : h \in \Classhp \big\}.
$$

\brem
%The set $\Classmp $ represents all bounded $\gg$-martingales  $Y$ for which  $Y_\tau$ is equal to the value at time $\tau$ of a $\gg$-predictable process.
(a) Due to the equality $\G_{\tau-}=\F_{\tau-}$, the set $\Classmp $ represents all bounded $\gg$-martingales  $Y$ for which  $Y_\tau$ is equal to the value at time $\tau$ of a $\gg$-predictable process, i.e., $Y_\tau \in \G_{\tau-}$.\\
(b) Since  the equality $\G_\tau= \F_\tau$ is not valid, in general,   $\Classm$ may differ from the class of all bounded $\gg$-martingales stopped at $\tau $. Note that $\Classm$ is equal
to the class of all bounded $\gg$-martingales stopped at $\tau $  if $\tau $ is an $\ff$-stopping time.
If $\G_\tau= \G_{\tau-}=\F_\tau$, then we have equality between the two classes $\Classmp$ and $\Classm$.   In general, one has only $\G_{\tau-}\subset \G_{\tau}$  where the inclusion may be strict. For instance, if $\tau$
is the first jump  of  a compound Poisson process $X_t=\sum _{n=1}^{N_t} Y_n$, the random variable $Y_1$ belongs to
$\G^X_\tau$ but not to $\G^X_{\tau-}$ where $\gg^X$ is the natural filtration of $X$.\\
%The set $\Classm$ is included in the class of all bounded $\gg$-martingales stopped at $\tau $.
 (c) From the credit risk perspective the computation of the quantity $\E(h_\tau \vert \G_t)$ is related to pricing of a recovery $h_\tau$ paid at time $\tau$. The representation property and the explicit form of the integrators allow us to find the hedging strategy (associated with hedging instruments obtained through the basic martingales): one needs either two or three traded assets, depending on the multiplicity of the filtration.
 \erem

We address the following general question: under which assumptions on the filtration $\ff$ and the random time $\tau $, any $\gg$-martingale $Y^h \in \Classm$ (or $Y^h \in \Classmp \subset \Classm $) admits an integral representation with respect to some \emph{fundamental} $\gg$-martingales? Of course, an essential part of this problem is a judicious specification of the family of fundamental $\gg$-martingales, which will serve as integrators in representation results for $\gg$-martingales. To be more specific, we wish to characterize the dynamics of any process  $Y^h \in \Classm $ (or $Y^h \in \Classmp $) and to obtain sufficient conditions for some integral representations to hold in the  filtration $\gg$ with respect to the $\gg$-martingale $\mt $ given by \eqref{mtau} and some auxiliary $\gg$-local martingales where the choice of these auxiliary martingales may depend on the setup at hand.

\bl \label{lemma1.1}
%(i)
If the process $h$ belongs to the class $\Classhp $ or if the process $h$ belongs to the class $\Classh $ and $\ao =\ap $, then the $\gg$-martingale
$Y^h_t :=  \E(h_\tau \vert \G_t)$ for $t \in \rr_+$ satisfies
\begin{equation} \label{y}
Y_t^h =
\ind_{\{\tau \leq t\}} h_\tau +  \ind_{\{t<\tau \} }
Z^{-1}_t \, \E\bigg( \int_t^{\infty } h_s \, d \ap_s+ h_\infty Z_\infty \, \Big|
\, \F_t \bigg)= A_t h_\tau +  (1-A_t) Z^{-1}_t X_t^h
\end{equation}
where we denote
\begin{equation} \label{yx1}
X_t^h:=\E \bigg( \int_t^{\infty} h_s \, d \ap_s +h_\infty Z_\infty \, \Big| \, \F_t \bigg) = \mx_t - \int_0^t h_s \, d \ap_s
\end{equation}
and $\mx$ is the uniformly integrable $\ff$-martingale given by
\begin{equation}\label{yx2}
\mx_t := \E \bigg( \int_0^{\infty} h_s \, d \ap_s +h_\infty Z_\infty \, \Big| \, \F_t \bigg).
\end{equation}
\el
\proof
If the process $h$ belongs to the class $\Classhp $, equality \eqref{y} was established in Elliott et al. \cite{EJY}.
They first show in Lemma 3.1 in \cite{EJY} that for any process $h \in \Classh$ one has
(recall that $Z>0$ on $[\![ 0 , \tau [\![$)
\begin{equation}\label{yvv1}
Y_t^h = \ind_{\{\tau \leq t\}} h_\tau +  \ind_{\{t<\tau \} } Z^{-1}_t \, \E\big( h_{\tau } \ind_{\{t<\tau \} } \, | \, \F_t \big).
\end{equation}
Subsequently, using the definition of dual predictable projection, it follows that for any process $h \in \Classhp$
\begin{equation} \label{yvv2}
\E\big( h_{\tau } \ind_{\{t<\tau<\infty \} } \, | \, \F_t \big)=\E\bigg( \int_t^{\infty } h_s \, d \ap_s \, \Big| \, \F_t \bigg).
\end{equation}
By combining \eqref{yvv1} with \eqref{yvv2}, we obtain \eqref{y}.

If the process $h$ belongs to the class $\Classh $ and $\ao =\ap $, using first the properties of the dual optional projection and subsequently the equality  ${A}^{o}= A^{p}$, we obtain, for any process $h \in \Classh $,
\begin{equation} \label{yvv3}
\E\big( h_{\tau } \ind_{\{t<\tau<\infty \} } \, | \, \F_t \big)=\E\bigg(
\int_t^{\infty }h_s\,d{A}^{o}_s \, \Big| \, \F_t \bigg) =
\E\bigg( \int_t^{\infty }h_s\,d\ap_s \, \Big| \, \F_t \bigg).
\end{equation}
By combining the last equality with \eqref{yvv1}, we conclude once again that \eqref{y} holds.
\endproof

\brem The condition $A^o=A^p$ is required for technical reasons. We have seen that, if either $\tau$ avoids $\ff$ stopping times or $\ff$ is a continuous filtration, then this condition holds. 
As derived in \cite{acj_thin}, the condition $A^o=A^p$ holds if and only if the random time $\tau$ satisfies that:
\begin{enumerate}
\item %avoids all $\ff$-totally inaccessible stopping times, i.e., 
$\P(\tau=T<\infty)=0$ for any $\ff$-totally inaccessible stopping time $T$,
\item $\P(\tau=S<\infty|\F_S)=\P(\tau=S<\infty|\F_{S-})$ for any $\ff$-predictable stopping time $S$.
\end{enumerate}
In particular the condition (2) is always satisfied in quasi-left continuous filtrations since then $\F_S=\F_{S-}$ for $\ff$-predictable stopping time $S$.
Therefore, if $\ff$ is a Poisson filtration (which is quasi-left continuous), the condition $A^o=A^p$ holds if and only if $\P(\tau=T_n<\infty)=0$ for any $n$ since any $\ff$-totally inaccessible stopping time $T$ satisfies $[\![T]\!]\subset \bigcup_n [\![T_n]\!]$.
One can  find such an example in \cite[Proposition 4]{aksamit/choulli/deng/jeanblanc}.
\erem

\section{Enlargement of a Poisson Filtration with a Random Time} \label{sec3}
%%%%%%%%%%%%%%%%%%%%%%%%%%%%%%%%%%%%%%%%%%%%%%%%%%%%%%%%%%%%%%%%%%%%%%%%%%%%%%%%%
%%%%%%%%%%%%%%%%%%%%%%%%%%%%%%%%%%%%%%%%%%%%%%%%%%%%%%%%%%%%%%%%%%%%%%%%%%%%%%%%%

Let $N$ be a Poisson process on $(\Omega, \G, \P)$ with intensity $\lambda$ and sequence of jump times denoted by $(T_n)_{n=1}^{\infty }$. We take $\ff = \ff^N$ to be the filtration generated by the Poisson process $N$, completed, and we denote by $M_t := N_t - \lambda t $ the compensated Poisson process, which is known to be an $\ff$-martingale. We consider
 an arbitrary random time $\tau$. Recall that the equality $Z = \mu - \ap$
is the Doob-Meyer decomposition of the Az\'ema supermartingale $Z$ of $\tau $.
It is worth noting that the Az\'ema supermartingale of any random
time $\tau $ is necessarily a process of finite variation when the filtration $\ff$ is generated
by a Poisson process.

From the predictable representation property of the compensated Poisson process (see, for instance, Proposition 8.3.5.1 in \cite{JYC}), the $\ff$-martingales $\mz $, $\mx $ and $m$ admit the integral representations
\be \label{mmi}
\mz_t=1 +\int_0^t \phi_s\, dM_s , \quad \mx_t=\mx_0+\int_ 0^t \phi^h_s\, dM_s , \quad m_t=1+\int_ 0^t \gamma_s\, dM_s
\ee
for some $\ff$-predictable processes $\phi , \phi^h $ and $\gamma$.

%%%%%%%%%%%%%%%%%%%%%%%%%%%%%%%%%%%%%%%%%%%%%%%%%%%%%%%%%%%%%%%%%%%%%%%%%%%%%%%%%%%%%
\subsection{General Predictable Representation Formula}  \lab{sec3.1}
%%%%%%%%%%%%%%%%%%%%%%%%%%%%%%%%%%%%%%%%%%%%%%%%%%%%%%%%%%%%%%%%%%%%%%%%%%%%%%%%%%%%%

The main difficulty in establishing the PRP for a
progressive enlargement is due to the fact that the jumps of
$\ff$-martingales may overlap with the jump of the process $A$. It
appears that even when the filtration $\ff$ is generated by a
Poisson process $N$, the validity of the PRP for the progressive
enlargement of $\ff$ with a random time $\tau $ is a challenging
problem for the part after $\tau$ if no additional assumptions are
made. Indeed, since the inclusion $\{ \Delta A >0 \} \subset \{ \Delta N >0\}$
may fail to hold, in general, it is rather hard to control a possible
overlap of jumps of processes $N$ and $A$ when the only
information about the $\ff$-conditional distribution of the random
time $\tau $ is the knowledge of its Az\'ema supermartingale $Z$.

Nevertheless, in the main result of this section (Theorem \ref{pox2}) we offer a
general representation formula for a $\gg$-martingale from $\Classmp$.
%It is fair to acknowledge that we need to introduce for this purpose an additional martingale to compensate for a potential mismatch between the jumps of the processes $A$ and $N$.
Subsequently, we
illustrate our general result by considering some special cases.

From Proposition \ref{proq}, we deduce that the compensated
martingale $M_t := N_t - \lambda t$ stopped at $\tau $ admits the following semimartingale decomposition with respect to $\gg$
\be \label{j1x}
\widehat M_t= M_{t \wedge \tau}-\int _0^{t\wedge \tau}\frac{d \cro {M,m}^\ff_s}{Z_{s-}}
= M_{t \wedge \tau}- \lambda \int _0^{t\wedge \tau} \frac{\gamma_s}{Z_{s-}}\, ds
\ee
where $\widehat M$ is a $\gg$-local martingale.
Note that the compensated Poisson process $M_t = N_t - \lambda t$ is an $\ff$-adapted (hence also $\gg$-adapted) process of finite variation so that it is manifestly a $\gg$-semimartingale for any choice of the random time $\tau$. Consequently, due to the PRP of the compensated Poisson process with respect to its natural filtration, no additional assumptions
regarding the random time $\tau$ are needed to ensure that the hypothesis $(H^\prime)$ is satisfied, that is,
any $\ff$-martingale is a $\gg$-semimartingale.

We define $A^*_t :=  \int_0^t \Delta N_s \, dA_s$ and
$\sigma :=\inf\, \{t \in \rr_+ :A^*_t=1\}$ so that
$A^*_t=\ind_{\{\sigma \leq t\} }$ and the jump times of the
process $A-A^*$ are disjoint from the sequence $(T_n)_{n=1}^{\infty}$ of jump times of $N$. The
$\gg$-adapted and increasing process $A^*$ is stopped at $\tau $ and it admits
a $\gg$-predictable compensator. Therefore, there exists a unique
$\ff$-predictable, increasing process $\Lambda^*$ such that the
process $ M^*_t:=A^*_t - \Lambda ^*_{t\wedge \tau} $ is a $\gg$-martingale stopped at $\tau $.
The $\gg$-compensator $\Lambda ^*$ is provided in the next lemma.
In the particular case where the graph of the random time $\tau$ is included in $\cup_{n=1}^\infty [\![T_n]\!]$, then one simply has that $M^*=\mt$.

\bl
\label{lem:Mxi}
Let  $A^*_t :=  \int_0^t \Delta N_s \, dA_s$, then
\be
\label{Mstar}
M^*_t=A^*_t -\int_0^t \frac{1-A_{s-}}{Z_{s-}}\, (\gamma_s-\phi_s)\, ds
\ee
is a $\gg$-local martingale. Let $\xi$ be an $\ff$-predictable process. Then the process $\wt{M}^\xi$ defined as
\be \lab{mmccx} \wt{M}^\xi_t :=
\xi_\tau\Delta N_\tau A_t -\int_0^t \frac{1-A_{s-}}{Z_{s-}}\, \xi_s(\gamma_s-\phi_s)\, ds
=\int_0^t \xi_s\, dM^*_s
\ee
is a $\gg$-local martingale.
\el

\proof
The assertion follows by Corollary \ref{cor1}, Lemma \ref{compensator} and \eqref{mmi}.
\endproof

\begin{Remark} Note that if $\tau$ avoids all $\ff$-stopping times, then $A^*=0$ and, since $m=\mu$, one has $\gamma=\phi$. Furthermore, using Corollary \ref{cor1} for $X=M$, we see that $\int_0^{\cdot \wedge \tau} \frac{1}{Z_{s-}}(\gamma_s-\phi_s) \lambda \, ds  =  \int _0^{\cdot \wedge \tau}  \frac{1}{Z_{s-}} \, dJ_s$ where $J$ is an increasing process, being the dual predictable projection of the increasing process $A\Delta M_\tau$. Therefore, the process $\gamma- \phi$ is nonnegative, as it must be. \end{Remark}

We are in a position to prove the main result of the paper. It should be stressed that the $\gg$-martingales $\mt$, $\wh{M}$ and $M^*$  (as given by equations \eqref{mtau}, \eqref{j1x} and \eqref{Mstar}, respectively) in the statement of Theorem \ref{pox2} are universal, in the sense that they depend on $N$ and $\tau$, but they are independent of the choice of the process $h \in \Classhp$. The following result furnishes the predictable representation for any $\gg$-martingale $Y^h$ from the class $\Classmp $.

\bt \lab{pox2}
If $\ff$ is a Poisson filtration, then any $\gg$-martingale $Y^h \in \Classmp $ admits the following predictable representation
\begin{align} \label{final}
dY^h_{t}&=\kappa^1_t\,d\mt_{t}+\kappa^2_t\,d\widehat M_{t}+\kappa^3_t\, dM^*_{t}
\end{align}
where
\begin{align*}
\kappa^1_t&=\left(h_{t}-\frac
{X^h_{t-}} {Z_{t-}}\right)
\left (1+\frac{\Delta \ap _{t}}{Z_{t-}-\Delta \ap_{t}}\ind_{\{Z_{t-}>\Delta \ap _{t}\}} \right)+\frac{X^h_{t-}}{Z_{t-}}\,\ind_{\{Z_{t-}=\Delta \ap_t>0\}}\\
\kappa^2_t&=\frac{1-A_{t-}}{Z_-} \left[ \left( \phi_{t}^h -\phi_{t} \frac{X^h_{t-}}
{Z_{t-}}\right)(1-\psi_t\phi_t)+\phi_t \frac{X^h_{t-}}{Z_{t-}} \, \ind_{\{Z_{t-}+\phi_t=0\}} \right] \\
\kappa^3_t&=-\psi_t\left (\phi^h_t-\phi_t\frac{X^h_{t-}}{Z_{t-}}\right)+\frac{X^h_{t-}}{Z_{t-}}\,\ind_{\{Z_{t-}+\phi_t=0\}}
\end{align*}
where we denote $\psi_t=\frac{1}{Z_{t-}+\phi_{t}}\,\ind_{\{Z_{t-}+\phi_{t}\neq 0\}}.$
\et

\proof Since only $\gg$-martingales stopped at $\tau$ are considered, it is enough to establish the decomposition on $[\![0,\tau]\!]$ and thus we work on $[\![0,\tau]\!]$ in the proof. For brevity, the variable $t$ is dropped.

\noindent{\it Step 1.}$\,$ In the first step, we derive representation \eqref{yeq3a}. From \eqref{y}, we obtain $Y^h = A h_\tau +  (1-A) V$ where $V := X^h Z^{-1}\ind_{\{Z>0\}}$. The It\^o formula yields
\begin{align} \label{popo}
dY^h &= h\, d\myh +(1-\myh_{-})\, dV
- V_{-} \, d\myh - \Delta V \Delta \myh \nonumber \\
&=(h- V_{-}) \, d\myh +(1-\myh_{-})\,
dV - \Delta V \Delta \myh
\end{align}
where by Lemma \ref{l:xz}(b)
\begin{align*} % \label{xpopo1}
dV& = \frac{1}{Z_{-}}\,dX^h
-\frac{V_{-}}{Z_{-}}\,dZ-\Delta V\,\frac{\Delta Z}{Z_{-}}.
\end{align*}
Recalling that $dX^h = d\mx - h \, d\ap$ and $dZ = d\mu- d\ap$, we obtain after elementary computations
\begin{equation}
\label{yeq1}
\begin{aligned}
dY^h&=(h-V_{-})\, d\mt +\frac{1-A_{-}}{Z_{-}}(d\mx-V_{-}\, d\mz)- (1-A_{-})\Delta V\,\frac{\Delta Z}{Z_{-}}-\Delta V\Delta A \\&=(h-V_{-})\,d\mt+\frac{1-A_{-}}{Z_{-}}\,\alpha\,dM-(1-A_{-})\Delta V\,\frac{\Delta Z}{Z_{-}}-\Delta V \Delta A
\end{aligned}
\end{equation}
where $\alpha=\phi^h-V_{-}\phi $. Observe that $\Delta Z=\Delta \mz-\Delta \ap $ and $\Delta X^h=\Delta \mx-h \Delta \ap $.
In view  of \eqref{mmi}, it is clear that $\{ \Delta \mz \ne 0 \} \subset \{\Delta N >0\}$ and $\{\Delta \mx \ne 0\} \subset \{ \Delta N > 0 \}$. Since the process $\ap $ is $\ff$-predictable and the jump times of $N$ are $\ff$-totally inaccessible, we have that $\Delta N \Delta \ap=0$. Using again \eqref{mmi}, we thus obtain
\bde
\Delta \mz \Delta \ap = \Delta \mx \Delta \ap=0.
\ede
%Using again \eqref{mmi}, we obtain
%\bde
%d[\mz ,\mx ] - d\cro{\mx , \mz}^\ff = \phi^h \phi \, dM , \quad
%d[\mz ,\mz ] - d\cro{\mz , \mz}^\ff = \phi^2 \, dM ,
%\ede
%since, in view of \eqref{mmi}, we have  $d(\Delta \mz \Delta \mx ) = d[\mz ,\mx ]$ and
%$d(\Delta \mz )^2 = d[\mz ,\mz ]$.
Hence
\begin{align*}
\Delta V&=\ind_{\{Z>0\}}\frac{1}{Z}\left ( \Delta \mx\  - V_{-}\Delta \mu-\left(h-V_-\right)\Delta \ap\right) -\ind_{\{Z=0\}}V_-\\
&= \psi (\phi^h - V_{-}\phi ) \Delta N
- \left( h- V_{-}\right)\frac{\Delta \ap }{(Z_{-}-\Delta \ap )}\ind_{\{Z_{-}>\Delta \ap \}}-\ind_{\{Z=0\}}V_-\\
&= \psi\alpha \Delta N - ( h- V_{-})\beta -\ind_{\{Z=0\}}V_-
\end{align*}
where
\begin{align*}
\beta  = \frac{\Delta \ap }{Z_{-}-\Delta \ap}\ind_{\{Z_{-}>\Delta \ap \}}
\end{align*}
so that
\begin{align} \label{yeq2}
\Delta V\, \frac{\Delta Z}{Z_{-}}  =&\psi \phi \alpha
  \frac{dN}{Z_{-}}+ \left( h-
V_{-}\right)\beta\,\frac {d\ap}{Z_{-}}-\ind_{\{Z=0\}}V_-\,
\frac{\Delta Z}{Z_{-}}
\end{align}
and
\begin{align} \label{yeq3}
\Delta V \Delta A = \psi \alpha \Delta N\,dA - (h-V_{-})\beta\,dA - \ind_{\{Z=0\}}V_-\,dA.
\end{align}
By combining \eqref{yeq1}, \eqref{yeq2} and \eqref{yeq3},  we obtain
\begin{align*}
dY^h& =(h-V_{-})\,d\mt+\frac{1-A_{-}}{Z_{-}}\,\alpha\,dM-(1-A_{-})\Delta V\,\frac{\Delta Z}{Z_{-}}-\Delta V\Delta A\\
& =(h-V_{-})\, d\mt + \frac{1-A_{-}}{Z_{-}}\,\alpha\,dM - \frac{1-A_{-}}{Z_{-}}\psi\phi\alpha\,dN \\
& - (1-A_{-})(h- V_{-})\beta\,\frac {d\ap}{Z_{-}}+(1-A_{-}) \ind_{\{Z=0\}}V_-\,
\frac{\Delta Z}{Z_{-}}- \psi\alpha\Delta N\, dA \\
&+(h-V_{-})\beta \, dA + \ind_{\{Z=0\}}V_-\,dA .
\end{align*}
Consequently, using the definitions of $\mt$ and $M$, we get
\begin{equation}
\label{yeq3a}
\begin{aligned}
dY^h&=(h-V_{-})(1+\beta )\,d\mt+\frac{1-A_{-}}{Z_{-}}\,\alpha (1-\psi\phi)\,dN-\frac{1-A_{-}}{Z_{-}}\psi\phi\alpha\lambda\, dt \\ & + (1-A_{-}) \ind_{\{Z=0\}}V_-\,\frac{\Delta Z}{Z_{-}}-\psi\alpha\Delta N\,dA +\ind_{\{Z=0\}}V_-\,dA.
\end{aligned}
\end{equation}

\noindent{\it Step 2.}$\,$ In this step, we analyze the term with $\Delta Z$ in  \eqref{yeq3a} and we establish representation  \eqref{yeq5}. Recall that $Z=Z_0+\int_0^{\cdot} \phi\,dM - \ap$ and $\Delta M = \Delta N$ so that $\Delta M \Delta \ap=0$.
Hence if $\Delta \ap>0$, then $\Delta Z = - \Delta \ap$ and thus
\begin{align*}
\{Z_{-}=\Delta \ap>0\}= \{Z_{-}-\Delta \ap=0,\, Z_{-}=\Delta \ap>0 \}=\{Z=0,\, Z_{-}=\Delta \ap>0\}.
\end{align*}
Therefore,
\begin{align*}
(1-A_{-}) \ind_{\{Z=0\}}V_-\,
\frac{\Delta Z}{Z_{-}} &= (1-A_{-}) \ind_{\{Z=0,\,Z_{-}=-\phi>0\}} \phi V_-\,
\frac{\Delta M}{Z_{-}} - (1-A_{-}) \ind_{\{Z=0,\,Z_{-}=\Delta \ap>0\}} V_-\,
\frac{\Delta \ap}{Z_{-}}\\
% &= (1-A_{-}) \ind_{\{Z_{-}=-\phi>0\}} \phi V_-\,
%\frac{\Delta N}{Z_{-}} - (1-A_{-}) \ind_{\{Z=0,\,Z_{-}=\Delta \ap>0\}} V_-\,
%\frac{\Delta \ap}{Z_{-}}\\
 &= (1-A_{-}) \ind_{\{Z_{-}=-\phi\}} \phi V_-\,
\frac{\Delta N}{Z_{-}} - (1-A_{-}) \ind_{\{Z_{-}=\Delta \ap>0\}} V_-\,
\frac{\Delta \ap}{Z_{-}}.
\end{align*}
We have thus shown that
\begin{align} \label{yeq4}
(1-A_{-}) \ind_{\{Z=0\}}V_-\,
\frac{\Delta Z}{Z_{-}}= (1-A_{-}) \phi  \myeta V_-\, \frac{dN}{Z_{-}} - (1-A_{-}) \myzeta V_-\, \frac{d\ap}{Z_{-}}
\end{align}
where we denote $\myeta :=\ind_{\{Z_-+\phi=0\}}$ and $\myzeta := \ind_{\{Z_{-}=\Delta \ap>0\}} =
\ind_{\{Z=0,\, Z_{-}=\Delta \ap>0\}}$.
It is obvious that the last term in \eqref{yeq3a} can be represented as follows
\begin{align*}
 \ind_{\{Z=0\}}V_-\,dA=( \ind_{\{Z=0\}}-\myzeta ) V_-\,dA+ \myzeta  V_-\,dA
\end{align*}
and thus, from \eqref{yeq3a} and \eqref{yeq4}, we obtain
\begin{equation}
\label{yeq5}
\begin{aligned}
dY^h& =\big[ \left(h-V_{-}\right) (1+\beta )+ \myzeta V_- \big] \, d\mt
+\frac{1-A_{-}}{Z_{-}}\, \big[ \alpha (1 -\psi \phi) + \phi  \myeta V_{-}\big] \, dN
\\&- \frac{1-A_{-}}{Z_{-}} \psi \phi \alpha \lambda \, dt - \psi \alpha \Delta N \, dA +( \ind_{\{Z=0\}}-\myzeta) V_-\,dA .
\end{aligned}
\end{equation}

\noindent{\it Step 3.}$\,$
To complete the derivation of \eqref{final}, it suffices to analyze the last term in \eqref{yeq5}.
Recall that $Z$ is a nonnegative, bounded, $\ff$-supermartingale and thus $Z=0$ on the set $\{Z_{-}=0\}$.
Hence, using the decomposition $Z=\mu - \ap$, we obtain % for every $t \in \rr_+$
\begin{align*}
\E \Big( \int_0^\infty \ind_{\{Z_{-}=Z=0\}} \, dA \Big)=
\E \Big( \int_0^\infty \ind_{\{Z_{-}=0\}}  \, dA \Big)=
\E \Big( \int_0^\infty \ind_{\{Z_{-}=0\}} \, d\ap \Big)=\E \Big( \int_0^\infty \ind_{\{Z_{-}=0\}} \, dZ \Big)=0
\end{align*}
%\begin{align*}
%\E \Big( \int_0^\infty \ind_{\{Z_{-}=Z=0\}} \, dA \Big)&=
%\E \Big( \int_0^\infty \ind_{\{Z_{-}=0\}}  \, dA \Big)=
%\E \Big( \int_0^\infty \ind_{\{Z_{-}=0\}} \, d\ap \Big)\\
%&=\E \Big( \int_0^\infty \ind_{\{Z_{-}=0\}} \, dZ \Big)=0
%\end{align*}
so that the nonnegative measure $\ind_{\{Z_{-}=Z=0\}}\,dA$ on $(\rr_+, \mathcal{B}(\rr_+))$ is null, almost surely.
% on the $\ff$-predictable $\sigma$-field.
Consequently,
\begin{align*} % \label{yeq6}
&\big( \ind_{\{Z =0\}} -\ind_{\{Z=0,\,Z_{-}=\Delta \ap>0\}}  \big) V_- \, dA
=\big( \ind_{\{Z=0,\, Z_{-}=-\phi>0\}} + \ind_{\{Z_{-}=Z=0 \}}  \big) V_- \, dA   \\
&=\ind_{\{Z=0,\, Z_{-}=-\phi>0\}} V_- \, dA
= \ind_{\{\Delta N=1,\, Z_{-}+\phi =0\}} V_- \Delta N \, dA = \ind_{\{Z_{-}+\phi =0\}} V_- \Delta N \, dA
= \myeta V_- \Delta N \, dA .\nonumber
\end{align*}
By substituting the last equality into \eqref{yeq5}, we obtain
\begin{align*}
dY^h& =\big[ \left(h-V_{-}\right) (1+\beta )+ \myzeta V_- \big] \, d\mt
+\frac{1-A_{-}}{Z_{-}}\, \big[ \alpha (1 -\psi \phi) + \phi  \myeta V_{-}\big] \, dN \\
&- \frac{1-A_{-}}{Z_{-}} \psi \phi \alpha \lambda \, dt +( -\psi \alpha + \myeta V_{-}) \Delta N \, dA .
\end{align*}
Using \eqref{j1x} and applying \eqref{mmccx} to the $\ff$-predictable process $\xi=-\psi \alpha + \myeta V_{-}$, we get
\begin{equation}
\label{final1}
\begin{aligned}
dY^h& =\big[ \left(h-V_{-}\right) (1+\beta )+ \myzeta V_- \big] \, d\mt
+\frac{1-A_{-}}{Z_{-}}\, \big[ \alpha (1 -\psi \phi) + \phi  \myeta V_{-}\big] \, d\wh{M}  \\
&+[ -\psi \alpha + \myeta V_{-}] \, dM^* + \frac{1-A_{-}}{Z_{-}}\, \lambda \, dK
\end{aligned}
\end{equation}
where the $\ff$-predictable process $K$ satisfies $K_0=0$ and
\begin{align*}
dK&= \big[ \alpha (1 -\psi \phi) + \phi  \myeta V_{-}\big] \big(1+ \frac{\gamma}{Z_{-}}\big) \, dt
- \psi \phi \alpha \, dt +( -\psi \alpha + \myeta V_{-})(\gamma-\phi )\, dt.
\end{align*}
%\begin{align*}
%dJ&= \big[ \alpha (1 -\psi \phi) + \phi  \myeta V_{-}\big] \big(\lambda + \frac{\lambda \gamma}{Z_{-}}\big) \, dt
%- \psi \phi \alpha \lambda \, dt +( -\psi \alpha + \myeta V_{-})(\gamma-\phi ) \lambda \, dt
%\end{align*}
The uniqueness of the Doob-Meyer decomposition implies that $K=0$ (this can also be checked by
direct computations). Since $K=0$, equation \eqref{final1} reduces to the desired representation  \eqref{final}.
\endproof

\brem \label{rem:taumi}
Let us comment on the relationship between Theorem \ref{pox2} and results obtained by Jeanblanc \& Song \cite{js3}.
Theorem 3.2 in \cite{js3} shows that the PRP holds for a class of $\gg$-local martingales stopped at $\tau$ with respect to $\mt$ and $\widehat M$ under the assumption that $\G_\tau\cap\{0<\tau<\infty\}=\G_{\tau-}\cap\{0<\tau<\infty\}$.
First, the latter condition is not always satisfied by a progressive enlargement of a Poisson filtration.
For example, let us take $\tau=T_1\ind_B+0.5 T_1\ind_{B^c}$ where $B\in\F_\infty$ but $B\notin\F_1$. Then $\G_{\tau-}=\sigma(\tau)$ and $\G_\tau=\sigma(\tau,B)$. Therefore, $\G_{\tau-}$ is strictly smaller than $\G_{\tau}$ and one cannot use results from \cite{js3}. To compute $\G_{\tau-}$ and $\G_\tau$, one can notice that $\gg$ is a marked point process
filtration and $\tau$ is a $\gg$ stopping time so that Theorems 2.2.14 and 2.2.15 in \cite{last} apply.
Second, in \cite{js3} the focus is on showing the existence of representation, whereas here we focus on explicit representations.
For this reason, we limit ourselves to a subclass of bounded $\gg$-martingales.
\erem

As a sanity check for representation \eqref{final} established in Theorem  \ref{pox2}, we will compute in two different ways the jump of $Y^h$ at $\tau$. The following remark will be used in the proof of Corollary \ref{gccd}.

\begin{Remark} \label{remarkxz0}
Let us define $R:=\inf\, \{t \in \rr_+ : Z_t=0 \}$. It is known that $X^h_R=0$ since $\tau\leq R$ and $X^h_R=\E(h_\tau \ind_{\{\tau>R\}} \vert \F_R)=0$ (see the proof of part (b) in Lemma \ref{l:xz}). This implies, in particular, that the equality $X^h_{\tau } = 0$ is valid on the event $\{Z_{\tau}=0\}$ since the equality $\tau = R$ holds on $\{Z_{\tau}=0\}$.
\end{Remark}

\bcor \label{gccd}
Under the assumptions of Theorem \ref{pox2}, let us set  $J_1 =\kappa^1_\tau \Delta \mt_\tau $
and $J_2=(\kappa^2_\tau +\kappa^3_\tau )\Delta N_\tau $ so that the jump  of $Y^h$ at time $\tau$ equals $J_1+J_2$. Then
\begin{align} \label{tyty}
J_1 = h_{\tau} -\frac {X^h_{\tau-}}{Z_{\tau-}} = \Delta Y^h_{\tau }, \quad J_2=0.
\end{align}
\ecor

\proof
On the one hand, in view of Lemma \ref{lemma1.1}, the jump of $Y^h$ at time $\tau$ satisfies $\Delta Y^h_{\tau }=h_\tau-\frac{X^h_{\tau-}}{Z_{\tau-}}$. On the other hand, from equation \eqref{final}, the jump  of $Y^h$ at time $\tau$  is
$$
\kappa^1_\tau \Delta \mt_\tau+\kappa^2_\tau\Delta \widehat M_\tau+\kappa^3_\tau \Delta M^*_\tau =
\kappa^1_\tau \Delta \mt_\tau + (\kappa^2_\tau +\kappa^3_\tau )\Delta N_\tau = J_1+J_2
$$
since $\Delta \widehat M_\tau=\Delta M_\tau=\Delta N_\tau$ and $\Delta M^*_\tau =\Delta N_\tau \Delta A_\tau=\Delta N_\tau$.
To compute $J_1$, we denote $I_1 =\ind_{\{Z_{\tau-}=\Delta \ap_{\tau} >0\}}$ and $I_2=\ind_{\{Z_{\tau-}>\Delta \ap_{\tau} \}} = \ind_{\{Z_{\tau-} \ne \Delta \ap_{\tau} \}}$ so that $I_1+I_2=1$.
Since $\Delta A_{\tau }=1$, it is obvious that
$$
I_1 \Delta \mt_\tau = I_1 \Big( 1 -\frac{\Delta \ap_{\tau }}{Z_{\tau-}}\Big)=0
$$
and thus
\begin{align} \label{ftft}
J_1 &= \kappa^1_{\tau} \Delta M^{(\tau)}_\tau = \bigg[ \Big(h_{\tau } -\frac{X^h_{\tau-}}{Z_{\tau-}} \Big) \Big( I_1+I_2+ \frac{\Delta \ap_{\tau}}{Z_{\tau-} -\Delta \ap_{\tau }}\, I_2 \Big)+\frac{X^h_{\tau-}}{Z_{\tau-}} \, I_1 \bigg]
 \Big( 1 -\frac{\Delta \ap_{\tau} }{Z_{\tau-}} \Big) \nonumber \\
&= \Big(h_{\tau } -\frac{X^h_{\tau-}}{Z_{\tau-}} \Big) \Big(1+ \frac{\Delta \ap_{\tau}}{Z_{\tau-} -\Delta \ap_{\tau }}\, \Big) \Big( 1 -\frac{\Delta \ap_{\tau} }{Z_{\tau-}} \Big) I_2
 = \Big(h_{\tau} -\frac {X^h_{\tau-}}{Z_{\tau-}}\Big) I_2 \\
 &= \Big(h_{\tau} -\frac {X^h_{\tau-}}{Z_{\tau-}}\Big)
  - \Big(h_{\tau} -\frac {X^h_{\tau-}}{Z_{\tau-}}\Big) I_1. \nonumber
\end{align}
Recall that $\Delta N \Delta \ap=0$. Hence, on the event $\{Z_{\tau-}=\Delta \ap_{\tau}>0\}$, from $Z=\mu-\ap$ we obtain $Z_{\tau}=0$ and thus, from Remark \ref{remarkxz0}, $X^h_\tau=0$. Furthermore, from \eqref{yx1} and \eqref{yx2}, we get $X^h_\tau= X^h_{\tau-}- h_\tau \Delta \ap_\tau=0$ on $\{Z_{\tau-}=\Delta \ap_{\tau}>0\}$ so that
\begin{equation} \label{xx33}
\Big(h_{\tau}-\frac {X^h_{\tau-}}{Z_{\tau-}}\Big) I_1
= \Big(h_{\tau}-\frac {X^h_{\tau-}}{\Delta \ap_{\tau}}\Big)\, \ind_{\{Z_{\tau-}=\Delta \ap_{\tau} >0\}}=0.
\end{equation}
By combining \eqref{ftft} with \eqref{xx33},  we conclude that the first equality in \eqref{tyty} is valid.
%$$
%J_1 = h_{\tau} -\frac {X^h_{\tau-}}{Z_{\tau-}}.
%$$
To complete the proof, it now suffices to show that $J_2=0$. We denote $\widetilde{I}_1= \ind_{\{Z_{\tau-}+\phi_{\tau}=0\}}$ and $\widetilde{I}_2=\ind_{\{Z_{\tau-}+\phi_{\tau} \ne 0\}}$ so that $\widetilde{I}_1+\widetilde{I}_2=1$. Also, we write
$$
\eta^h_{\tau } = \Big(\phi^h_{\tau } -\phi_{\tau} \frac{X^h_{\tau-}}{Z_{\tau-}}\Big).
$$
On the event $\{Z_{\tau-}+\phi_{\tau}=0\}$ we have that $Z_{\tau-}=-\phi_{\tau}>0$ and thus $\frac{X^h_{\tau-}}{Z^2_{\tau-}} \phi_{\tau} = -\frac{X^h_{\tau-}}{Z_{\tau-}}$. Therefore,
\begin{align*}
J_2 &= (\kappa^2_\tau +\kappa^3_\tau )\Delta N_\tau = \bigg[ \frac{1}{Z_{\tau-}} \, \eta^h_{\tau} (1- \psi_{\tau} \phi_{\tau} ) + \frac{X^h_{\tau-}}{Z^2_{\tau-}}\, \phi_{\tau} \widetilde{I}_1 - \eta^h_{\tau}\psi_{\tau }
 + \frac{X^h_{\tau-}}{Z_{\tau-}} \, \widetilde{I}_1 \bigg] \Delta N_\tau \\
 &= \eta^h_{\tau} \bigg[ \frac{1}{Z_{\tau-}} \Big(\widetilde{I}_1+\widetilde{I}_2 -\frac{\phi_{\tau} }{Z_{\tau-}+\phi_{\tau}} \, \widetilde{I}_2 \Big) -  \frac{1}{Z_{\tau-}+\phi_{\tau}} \, \widetilde{I}_2 \bigg] \Delta N_\tau
=\frac{1}{Z_{\tau-}}\Big(\phi^h_{\tau } -\phi_{\tau} \frac{X^h_{\tau-}}{Z_{\tau-}}\Big)\widetilde{I}_1 \Delta N_\tau \\
& =  \frac{1}{Z_{\tau-}}\Big(\Delta X^h_{\tau } + X^h_{\tau-}\Big)\ind_{\{Z_{\tau-}+\phi_{\tau}=0\}} \Delta N_\tau
  = \frac{ X^h_{\tau }}{Z_{\tau-}}  \ind_{\{Z_{\tau}=0\}} \Delta N_\tau =0
\end{align*}
where we used \eqref{mmi} and where the last equality holds since, by virtue of Lemma \ref{l:xz} (see Remark \ref{remarkxz0}), we have that $X^h_{\tau } = 0$ on the event $\{Z_{\tau}=0\}$.
\endproof

Let us give an alternative proof for the second equality in \eqref{xx33}.
It is known (see, e.g., Corollary 2.8 in \cite{nk}) that if $U$ is a right-continuous $\gg$-local martingale and $\vartheta$ a $\gg$-predictable stopping time, then
$$\E( U_\vartheta \ind_{\{\vartheta <\infty\}}\vert \G_{\vartheta-}) = U_{\vartheta-}\ind_{\{\vartheta <\infty\}}.$$
Moreover, if $B$ is a $\gg$-predictable set and $D_B$ the d\'ebut of $B$, then $[D_B]\subset B$ implies that $D_B$ is a
$\gg$-predictable stopping time. Let us take $E= \{Z_- =\Delta \ap>0\}$. Then we see that $\tau^*:= D_E$ is a $\gg$-predictable stopping time and $\tau^* \ind_{\{\tau^*<\infty\}}=\tau \ind_{\{\tau^*<\infty\}}$. Since $Y^h_\tau=h_\tau$, we obtain
\begin{align*}
h_\tau  \ind_{\{\tau^*<\infty\}} &= h_{\tau^*}\ind_{\{\tau^*<\infty\}}=\E( h_{\tau*} \ind_{\{\tau^*<\infty\}} \vert \G_{\tau^*-})=\E( Y^h_{\tau^* }\ind_{\{\tau^*<\infty\}} \vert \G_{\tau^* -})\\
&=Y^h_{\tau^*-}  \ind_{\{\tau^*<\infty\}}=Y^h_{\tau -} \ind_{\{\tau^*<\infty\}}
\end{align*}
where we have also used the fact that $h_{\tau*} \ind_{\{\tau^*<\infty\}} \in \G_{\tau^*-}$ since $h$ is $\ff$-predictable.
We conclude that $h_\tau =Y^h_{\tau -}$ on $E$ and thus equality \eqref{xx33} holds meaning that $Y^h$ has no jump at $\tau$ on $E$.

The following result shows that the predictable representation established in Theorem \ref{pox2} simplifies under the additional
assumption that $Z>0$.

\bt \label{cor:posi}
Assume that $\ff$ is a Poisson filtration and $Z>0$. Then any $\gg$-martingale $Y^h \in \Classmp $ admits the following predictable representation
\begin{align}
\label{eq:mainth}
dY^h_t = \frac {Z_{t-}}{Z_{t-}-\Delta \ap_t} \left(h_t-\frac
{X^h_{t-}} {Z_{t-}}\right)d\mt_t +\frac {1-A_{t-}} { Z_{t-}+\phi_t
} \left( \phi^h _t-\phi_t \frac{X^h_{t-}}{Z_{t-}}\right) dM^\bot_t
\end{align}
where $M^\bot:= \widehat M - M^*$ is orthogonal to $\mt$.
%Here the $\ff$-predictable processes $\phi $ and $\phi^h$ are given by \eqref{mmi}, the $\gg$-martingale $M^*$  by \eqref{Mstar} and $\gamma=\ind_{\{Z_{-}>0\}}\ind_{\{Z_{-}+\phi>0\}}.$
\et

\proof
If $Z>0$, we have $\ind_{\{Z_{-} = \Delta \ap>0\}}=0$ and $\ind_{\{Z_-+\phi=0\}}\Delta N=0$. We have $d\widehat M_t= dN_t + \vartheta_t\, dt$ and $dM^*_t = \Delta A_t \, dN_t + \xi_t\, dt$ for some processes $\vartheta, \xi$ and it is now easy to deduce
from the equality $\ind_{\{Z_-+\phi=0\}}\Delta N=0$ that $ \ind_{\{Z_{t-}+\phi_t \ne 0\}} \vartheta_t\, dt=  \vartheta_t \, dt$ and $ \ind_{\{Z_{t-}+\phi_t \ne 0\}} \xi_t\, dt=  \xi_t\, dt$. Consequently, the equalities $dM^*_t = \ind_{\{Z_-+\phi \ne 0\}}\, dM^*_t
$ and $d\widehat M_t = \ind_{\{Z_-+\phi \ne 0\}}\, d\widehat M_t$ are valid, and thus a straightforward application of Theorem  \ref{pox2}  yields
\begin{align*}
\kappa^2_t\, d\widehat M_t +\kappa^3_t \, dM_t^*&=(1-A_{t-}) \left( \phi^h_t -\phi_t \, \frac{X^h_{t-}}{Z_{t-}}\right)
\Bigg(\frac{1}{Z_{t-}}\bigg(1- \frac{ \phi_t}{Z_{t-}+\phi_t} \bigg)\, d\widehat M_t- \frac{1}{Z_{t-}+\phi_t}\, dM^*_t \Bigg)
\\&= \left( \phi^h_t -\phi_t \, \frac{X^h_{t-}}{Z_{t-}}\right)\frac{1}{Z_{t-}+\phi_t} \, dM^\bot_t .
\end{align*}
The orthogonality of $M^\bot$ and $\mt$ follows from the fact that these purely discontinuous martingales do not have common jumps.
\endproof

\brem
The PRP is closely related to the notion of \emph{multiplicity}. Here by the multiplicity of the filtration $\gg$ for the class $\Classmp$, we mean the minimal number of mutually orthogonal martingales needed to represent all martingales from $\Classmp$ as stochastic integrals. Therefore, under the assumptions of Theorem \ref{cor:posi}, the multiplicity of $\gg$ for the class $\Classmp$ equals two. For more details on the multiplicity of a filtration, see Davis \& Varaiya \cite{dv}, Davis \cite{d}
and the references therein.
\erem

%\proof Since $\{ \Delta\mx \ne 0 \}
%\subset \{ \Delta N > 0 \}$ and $\{ \Delta\mz \ne 0 \} \subset \{
%\Delta N > 0 \}$, it is clear that
%$$
%\left (\Delta\mx_\tau-  \Delta\mz_\tau
%\frac{X^h_{\tau-}}{Z_{\tau-}}\right)dA_{} = \left (\Delta\mx_\tau-
%\Delta\mz_\tau\frac{X^h_{\tau-}}{Z_{\tau-}} \right)dA_{}^*
%=\left(\phi^h_{}-\phi_{}\frac{X^h_{-}}{Z_{-}}\right) dA_{}^*
%$$
%where the second equality follows from \eqref{j1x}. It is now easy to check that
%$$
%d\wt{M}^h_{} = dB^h_{}-\frac{1-A_{-}}{Z_{-}}\, d\Bhpff_{}= -\left
%(\phi^h_{}-\phi_{}\frac{X^h_{-}}{Z_{-}} \right)dM^*_{}
%$$
%and the assertion follows by Proposition \ref{pox2} below.
%\endproof

We present further consequences of Theorem \ref{pox2}. We first consider the case of independent random time.

\bcor \label{cor:one}
Assume that $\tau$ is independent of the filtration $\ff$ of the Poisson process $N$. Then for any process $Y^h \in \Classmp $ we have
\be \label{rep2}
dY^h_t =\frac {Z_{t-}}  {Z_{t-}-\Delta \ap_t} \left(h_t-\frac
{X^h_{t-}} {Z_{t-}}\right)\, d\mt_t +\frac {1-A_{t-}}{Z_{t-}} \phi_t^h \, dM_t
\ee
where the $\ff$-predictable process $\phi^h$ is given by \eqref{mmi}.
\ecor

\proof
We have $\Delta N \Delta A= 0$ so that $A^* = 0$ and $M^*= 0$. Obviously, $\ff$ is
immersed in $\gg$ and the Az\'ema supermartingale $Z$ is a decreasing deterministic function and thus it is
$\ff$-predictable. Hence $\mu =1$ so that $\phi=0$. An application of Theorem \ref{pox2} concludes the proof.
\endproof

In the next result, we postulate the avoidance property for  the random time $\tau$.

\bcor \label{corss}
Assume that $\ff$ is a Poisson filtration and the random time $\tau$ avoids all $\ff$-stopping times. Then for any process $Y^h \in \Classmp $ we have
\be \label{rep1}
dY^h_t = \left(h_t-\frac{X^h_{t-}}{Z_{t-}}\right)d\mt_t + \frac{1-A_{t-}}{Z_{t-}+\phi_t} \left( \phi_{t}^h -\phi _{t} \frac{X^h_{t-}}{Z_{t-}}\right) d\widehat M_t\,.
\ee
 \ecor

\proof
The postulated avoidance property implies that, on the one hand, $A^*=0$ and thus also $M^*=0$ and, on the other  hand, the process $\ap$  is continuous (see, e.g., Lemma 8.13 in Nikeghbali \cite{nk}). Hence equality \eqref{rep1} follows from  Theorem \ref{pox2}.
\endproof

%%%%%%%%%%%%%%%%%%%%%%%%%%%%%%%%%%%%%%%%%%%%%%%%%%%%%%%%%%%%%%%%%%%%
\subsection{Cox Construction Example} \lab{sec3.2}
%%%%%%%%%%%%%%%%%%%%%%%%%%%%%%%%%%%%%%%%%%%%%%%%%%%%%%%%%%%%%%%%%%%%

We study here the special case where the immersion property between $\ff$ and $\gg$ holds.
The filtered probability space $(\Omega, \G, \ff,\P)$ supports a random variable $\Theta$  with the unit exponential distribution and such that $\Theta $ is independent of the filtration $\ff$ generated by a Poisson process $N$. The random time $\tau$ is given by a {\it Cox construction},  more precisely,
\be \lab{bbg1}
\tau = \inf \, \{ t \in \rr_+ \! :\, \psi (N_t) \geq \Theta \}
\ee
for some non-decreasing function $\psi : \rr_+ \to \rr_+$ such that $\psi (0)=0$ and $\lim_{\, x \to \infty } \psi (x)= \infty .$

The Az\'ema supermartingale equals $Z_t =e^{- \psi (N_t)}$ for all $t \in \rr_+$ and thus it is decreasing,
but not $\ff$-predictable.   Since  the construction (\ref{bbg1}) implies that the filtration $\ff$ is immersed in $\gg$ (see Kusuoka \cite{KS} or Elliott et al. \cite{EJY}),
the compensated Poisson $\ff$-martingale $M$ is also a $\gg$-martingale. It is
crucial to observe that the inclusion $[\![\tau]\!] \subset
\cup_n[\![T_n]\!] $ holds, meaning that a jump of the process $A$
may only occur when the Poisson process $N$ has a jump, that is,
$\{ \Delta A >0 \} \subset \{ \Delta N >0\}$. We first compute
explicitly the Doob-Meyer decomposition of $Z$ and we show that $\ap $ is continuous.
\bl \lab{lm44} Let $Z = e^{- \psi (N)}$ where $N$ is a Poisson process.
Then the following assertions hold: \hfill \break (i) $Z$ admits
the Doob-Meyer decomposition  $Z=\mz -\ap$ where
\be \label{mmiy}
\mu_t= 1+ \int_0^t \big( e^{-\psi (N_{s-}+1)}-e^{-\psi (N_{s-})} \big) \, dM_s ,\quad
\ap_t = \int_0^t \lambda \big( e^{-\psi (N_{s})}-e^{-\psi (N_{s}+1)} \big) \, ds ,
\ee
(ii) the process
\be \label{mrmiy}
\mt_t :=
A_t- \int_0^{t\wedge \tau } e^{\psi (N_s)}  \, d\ap_s  = A_t- \int_0^{t\wedge \tau }  \lambda  \big( 1- e^{\psi (N_s)-\psi (N_{s}+1)} \big) \, ds
\ee
is a $\gg$-martingale and thus the random time $\tau $ is a totally inaccessible $\gg$-stopping time.
\el

\proof
Equation \eqref{mmiy} is an immediate consequence of elementary computations
\bde
Z_t= Z_0 + \sum_{0< s \leq t } \big( e^{-\psi (N_{s})}-e^{-\psi (N_{s-})} \big)
 = 1 + \int_0^t \big( e^{-\psi (N_{s-}+1)}-e^{-\psi (N_{s-})} \big) \, dN_s
\ede
and the equality $M_t = N_t - \lambda t$.
For the second part,  we note that the $\gg$-martingale property of the process $\mt$ is
a consequence of \eqref{mtau}, whereas the fact that $\tau $ is a
totally inaccessible $\gg$-stopping time follows from the continuity of the
compensator in \eqref{mrmiy}.
\endproof

Recall that the process $X^h$ is defined by \eqref{yx1} and note
that, due to the continuity of $\ap$, we have
\be \lab{ut66} \{ \Delta \mt >0 \}
= \{ \Delta A >0 \} \subset \{ \Delta N >0\} = \{ \Delta M >0\}.
\ee

Our goal is to show that an explicit representation is available for some martingales
from the class $\Classmp $ within the present setup, that is, when the Az\'ema supermartingale $Z$ is not $\ff$-predictable.  For the corresponding result when the immersion hypothesis with an $\ff$-predictable Az\'ema supermartingale holds, see Example \ref{ex4.2}.

\bp \label{proc1}
Consider the $\gg$-martingale $Y^h \in \Classmp$ where the process $h \in \Classhp $
is given by $h_t=h(N_{t-})$ for some function $h : \mathbb N \to \rr $. Then $Y^h$ admits the following predictable representation
\begin{equation} \lab{axx33}
dY^h_t = \big( h_t - \wt h(N_{t-}+1) \big) \,d\mt_t
+(1-A_{t-})\big( \wt h(N_{t-}+1)-\wt h(N_{t-})\big)\,dM_{t}
\end{equation}
or, equivalently,
\begin{equation} \lab{cxx33}
dY^h_t =\big( h_t-  \wt h (N_{t-}) \big) \, d\mt_t +(1-A_{t-})
\big( \wt h (N_{t-}+1) - \wt h (N_{t-})\big) \, dK^{(\tau
)}_t
\end{equation}
where the processes $\mt$ and $K^{(\tau )}:= M -\mt$ are orthogonal $\gg$-martingales and $\widetilde h(x) :=  e^{\psi (x)}  \Phi(x)$ where the function  $\Phi : \mathbb N \to \rr$ is given by
\be \lab{nhy}
 \Phi(x) :=  \int_0^\infty  \lambda \E  \Big( h(N_s +x) \big(e^{-\psi (N_{s}+x)}-e^{-\psi (N_{s}+x+1)} \big) \Big) ds .
\ee
\ep

\proof
In view of equation \eqref{y} and Lemma \ref{lm44}, we have
\bde
Y^h_t=  A_t h(N_{\tau -}) + (1-A_t)  e^{\psi (N_t)} \, \E \left(  \int_t^\infty
 \lambda h(N_s) \big( e^{-\psi (N_{s})}-e^{-\psi (N_{s}+1)} \big)\, ds \,\Big|\, \F_t \right).
\ede
Using the time-homogeneity and the independence of increments of the Poisson process $N$, we obtain
\be \label{vcvc}
Y^h_t = A_t h(N_{\tau -}) +(1-A_t) e^{\psi (N_{t})}  \Phi(N_t)
=  A_t h(N_{\tau -}) +(1-A_t) \widetilde h(N_t)
\ee
where we denote $\widetilde h(x) :=  e^{\psi (x)}  \Phi(x)$.
It is useful to observe that $X^h_t = \Phi (N_t)$ and for every $x \in \mathbb N$
\be \lab{inhy}
\Phi (x+1) =  \Phi(x) - h(x) \big( e^{-\psi (x)}-e^{-\psi (x+1)} \big).
\ee
Equality \eqref{inhy} holds, since from the properties of the Poisson process we deduce that for any function $g : \mathbb N \to \rr $ we have
\begin{align*}
\Psi(x) &:=  \int_0^\infty  \lambda \E ( g(N_s+x) )\,  ds = \int_0^\infty  \lambda \E ( g(N_s +x+1)) \,  ds
+ \int_0^\infty  \lambda \E \big( g(x)\ind_{\{N_s=0\}} \big) \,  ds \\
&= \Psi (x+1) + g(x) \int_0^\infty \lambda e^{-\lambda s}\, ds = \Psi (x+1) + g(x).
\end{align*}
 %\mb{that, since  \begin{align*}\Phi(x)&=  \E \Big(\int_0^{T_1}  \lambda   \big( h(N_s +x) \big(e^{-\psi (N_{s}+x)}-e^{-\psi (N_{s}+x+1)} \big)ds  \,\Big|\, \F_t\Big) \\&+  \E \Big( \int_0^\infty  \lambda   h(N_s +x) \big(e^{-\psi (N_{s}+x)}-e^{-\psi (N_{s}+x+1)} \big)ds\, \Big|\, \F_t \Big) \end{align*} and $\E(T_1)=1/\lambda$,}
Using \eqref{yx1} and \eqref{inhy}, we obtain
\begin{align*}
\mx_t & = X^h_t + \int_0^t h_s \, d\ap_s =
\Phi (N_t) - \int_0^t \lambda h(N_s)\big( e^{- \psi (N_s+1)} - e^{- \psi (N_s)} \big) \, ds \\
&= \Phi (0) + \int_0^t \big( \Phi (N_{s-}+1) - \Phi (N_{s-}) \big) \, dN_s - \int_0^t \lambda h(N_{s-})\big( e^{- \psi (N_{s-}+1)}
    - e^{- \psi (N_{s-})} \big) \, ds \\
&= \Phi (0) + \int_0^t h(N_{s-}) \big( e^{-\psi (N_{s-}+1)}-e^{-\psi (N_{s-})} \big) \, dN_s  - \int_0^t \lambda h(N_{s-})\big( e^{- \psi (N_{s-}+1)} - e^{- \psi (N_{s-})} \big) \, ds \\
& = \Phi (0) + \int_0^t h(N_{s-}) \big( e^{-\psi (N_{s-}+1)}-e^{-\psi (N_{s-})} \big) \, dM_s .
\end{align*}
From \eqref{mmiy} and the just derived representation for $\mx$, we deduce that the $\ff$-predictable processes $\phi$ and $\phi^h$
appearing in \eqref{mmi} are given by
\bde
\phi_t = e^{- \psi (N_{t-}+1)} - e^{- \psi (N_{t-})},
\quad  \phi^h_t = h(N_{t-}) \big( e^{-\psi (N_{t-}+1)}-e^{-\psi (N_{t-})} \big).
\ede
By substituting these processes into \eqref{eq:mainth}
and noting that $\Delta \ap =0$ (see Lemma \ref{lm44}), we obtain
\begin{align}
 dY^h_t & =  \big( h_t- \wt h (N_{t-}) \big)\, d\mt_t
 + (1-A_{t-}) e^{\psi (N_{t-}+1)} \big( \phi_t^h -\phi_t \wt h(N_{t-}) \big)\, dK^{(\tau )}_t \nonumber \\
 & = \big( h_t-  \wt h (N_{t-}) \big) \, d\mt_t
+(1-A_{t-}) \big( \wt h (N_{t-}+1) - \wt h (N_{t-})\big) \, dK^{(\tau )}_t \label{ttgg}
\end{align}
where  $K^{(\tau )}:=M-\mt =\wh M-\mt$. To establish the second equality in \eqref{ttgg}, we used the following elementary equality
\bde
e^{\psi (x+1)} \left[ h(x) \big( e^{-\psi (x+1)}- e^{-\psi (x)}\big) - \big( e^{-\psi (x+1)}- e^{-\psi (x)}\big)
\wt h(x) \right] = \wt h(x+1) - \wt h(x),
\ede
which follows from \eqref{inhy} and the definition of the function $\wt h$.
Note that formula \eqref{ttgg} coincides with the predictable representation \eqref{cxx33}.
\endproof

\section{Enlargement of a General Filtration with a Pseudo-Stopping Time} \lab{sec4}
%%%%%%%%%%%%%%%%%%%%%%%%%%%%%%%%%%%%%%%%%%%%%%%%%%%%%%%%%%%%%%%%%%%%%%%%%%%%%%%%%%%%%%%%%%%%%%%%%
%%%%%%%%%%%%%%%%%%%%%%%%%%%%%%%%%%%%%%%%%%%%%%%%%%%%%%%%%%%%%%%%%%%%%%%%%%%%%%%%%%%%%%%%%%%%%%%%%

In this final section, we relax the assumption that the filtration $\ff$ is
generated by a Poisson process and we consider instead a general
right-continuous and $\P$-completed filtration $\ff$. As before,
we search for an integral representation of a $\gg$-martingale
$Y^h$ in terms of integrals with respect to the $\gg$-martingale
$\mt$ given by \eqref{mtau} and some  auxiliary $\gg$-martingales
associated with $h$ and $\tau $ (typically, the $\gg$-martingale parts of the $\ff$-martingales $\mx $ and $\mz$ stopped at $\tau $).
We henceforth work under the standing assumption that the random time $\tau$ is an $\ff$-pseudo-stopping time.
Recall from Nikeghbali \& Yor \cite{NY} that a random time $\tau $ is called an $\ff$-{\it pseudo-stopping time} if for any bounded $\ff$-martingale
$M$ we have $M_0 = \E (M_{\tau })$. Since $\tau $ is a $\gg$-stopping time,
it is clear that the immersion property between $\ff$ and $\gg$ is stronger than the assumption that $\tau $ is an $\ff$-pseudo-stopping time. For the following properties of pseudo-stopping times, we refer to Theorem 1 in
Nikeghbali \& Yor \cite{NY}, Theorem 9.20 in Nikeghbali \cite{nk} (see also Theorem 3 in Aksamit \& Li \cite{al} for a generalization of these results to the case of random time that is not necessarily finite). We emphasize that, if $\tau$ is an $\ff$-pseudo-stopping time, and $X$ is an $\ff$-martingale, then the stochastic integral $\int_0^{\cdot \wedge \tau}H_s\,dX_s$ is well defined for a $\gg$-predictable process $H$ satisfying suitable integrability conditions.

\bp \label{pseudo}
The following conditions are equivalent: \\
(i) a random time $\tau $ is an $\ff$-pseudo-stopping time,\\
(ii) the equality $m = 1 $ holds, \\
(iii) the Az\'ema supermartingale of $\tau $ satisfies $Z = 1 - \ao $, \\
(iv) for every $\ff$-local martingale $M$, the stopped process $\overline M := M^{\tau }$ is a $\gg$-local martingale.
\ep

From part (iv) in Proposition \ref{pseudo}, it follows that if $\tau $ is an $\ff$-pseudo-stopping time, then for
every $\ff$-local martingale $M$ the stopped process $\overline M := M^{\tau }$ is a $\gg$-local martingale where $M_t^{\tau } := M_{t\wedge \tau }$ for all $t \in \rr_+$. Accordingly, $\bmz$ (resp.,
$\bmx $) stands for the $\gg$-local martingale obtained by
stopping at $\tau $ the $\ff$-local martingale  $\mz$ (resp.,
$\mx$) meaning that $\bmz := \mz^{\tau}$ (resp., $\bmx :=
(\mx)^{\tau}$). Note that in the case of pseudo-stopping time the
$\ff$-martingale $\mu$ has finite variation since $\mu=1-\ao-\ap$.

Let us recall a particular result due to Coculescu et al. \cite{cjn}, which is closely related to our studies
in this section. Under the assumptions that:  (i) $\tau$ avoids $\ff$-stopping times, (ii) $\ff$ is immersed in $\gg$, and (iii) there exists an $\ff$-martingale $M$ that enjoys the PRP for the filtration  $\ff$,  Proposition 5.7 in \cite{cjn} shows
that the pair $(M, M^{(\tau)})$ has the PRP for the filtration $\gg$.

%%%%%%%%%%%%%%%%%%%%%%%%%%%%%%%%%%%%%%%%%%%%%%%%%%%%%%%%%%%%%
\subsection{Independent Random Time} \lab{sec4.1}
%%%%%%%%%%%%%%%%%%%%%%%%%%%%%%%%%%%%%%%%%%%%%%%%%%%%%%%%%%%%%

Let us first examine a very particular case where the random time
$\tau $ is independent of the filtration $\ff$ (hence it is an
$\ff$-pseudo-stopping time). Of course, this setup covers also the
case of the trivial filtration $\ff$. For a more detailed study of
the case of the trivial filtration $\ff$, we refer to Chou \&
Meyer~\cite{ChM}.
Note that in the next result we are dealing
with general $\gg$-martingales, and not only with
$\gg$-martingales stopped at $\tau $. Since we do not assume that
$\P(\tau = t)=0$ for any $t$, this simple result is not a
consequence of the one of Coculescu et al. \cite{cjn}.

\bp Assume that $\tau$ is independent of $\ff$ and $Z>0$. Let $M$ be any
$\ff$-martingale such that the PRP holds for $\ff$. Then the
progressive enlargement $\gg$ has PRP with respect to the pair
$(\mt , M)$.
\ep

\proof
Let $\xi$ be a bounded,  $\F_\infty$-measurable random variable and let $h$ be a
bounded Borel function. Due to the assumed independence of $\tau $ and $\ff$,
we have that $Z_t= \P (\tau > t )$ for all $t\in \rr_+$.
Let us set $m^\xi_t:=\E(\xi\vert \F_t)$. Then the $\gg$-martingale $Y_t := \E( \xi h(\tau)\vert \G_t)$ (note that $Y$ is not stopped at $\tau $) satisfies
\begin{align*}
Y_t = h(\tau) m^\xi_t A_t + (1-A_t)  Z^{-1}_t \, \E(\xi h(\tau)\ind_{\{\tau > t\}} \vert \F_t)
 =h(\tau) m^\xi_t A_t + (1-A_t) m^\xi_t \widehat h(t)
\end{align*}
where $\widehat h(t) := - Z^{-1}_t \int_t^\infty h(s)\, dZ_s$. The integration
by parts formula yields
$$
d((1-A_t)\widehat h(t))=(1-A_{t-}) h(t) Z^{-1}_{t-}\, dZ_t -
(1-A_{t-})\widehat h(t) Z^{-1}_{t-} \, dZ_t - \widehat h(t) \, dA_t .
$$
Consequently,
\begin{align*}
dY_t&= m^\xi_{t-} h(t)\, dA_t+ m^\xi_{t-} d((1-A_t) \widehat h(t)\big)+\big(  h(\tau ) A_t +(1-A_t) \widehat h(t) \big)\, dm^\xi_t\\
&= m^\xi_{t-}  h(t)\, dA_t + m^\xi_{t-} (1-A_{t-}) h(t) Z^{-1}_{t-}\, dZ_t -
 m^\xi_{t-} \big((1-A_{t-}) \widehat h(t) Z^{-1}_{t-} \,dZ_t + \widehat h(t)\, dA_t \big)
 \\& + \big( h(\tau ) A_t +(1-A_t)\widehat h(t) \big)\, dm^\xi_t\\
 &=m^\xi_{t-} (h(t)-\widehat h(t) ) \, dM^{(\tau)}_t +\big( h(\tau ) A_t +(1-A_t)\widehat h(t)\big) \, dm^\xi_t .
\end{align*}
It is clear that $m^\xi_{t-} (h(t)-\widehat h(t) ) \, dM^{(\tau)}_t$ is a $\gg$-martingale, since
the process $m^\xi_{t-} (h(t)-\widehat h(t) )$ is $\ff$-predictable (hence also $\gg$-predictable).
In contrast, the process $ h(\tau ) A_t +(1-A_t)\widehat h(t)$ is not $\gg$-predictable, however, due to
the independence of $\tau$ and $\ff$, one can check that any integral of the form
$\int_0^t k_s(\tau )\, dx_s$ is a $\gg$-martingale if for any $u \in \mathbb{R}_+$ the process $k(u)$ is bounded and
$\ff$-predictable and $x$ is an $\ff$-martingale. Choosing $k_s(u)=h(u)\ind_{\{u\leq s \}}+ \widehat h(s)\ind_{\{s<u\}}$ leads to the martingale property of  the integral with respect to $m^\xi$.
If $M$ is some $\ff$-martingale satisfying the PRP for the filtration $\ff$, then
$dm^\xi_t= \theta_t\, dM_t$ for some $\ff$-predictable process $\theta $. Moreover, since the immersion is satisfied,
$M$ is also a $\gg$-martingale. It is now clear that the pair $(M^{(\tau)},M)$ has indeed the representation property with respect to $\gg$.
\endproof

In general, the study of the problem of a representation after $\tau$ under the postulate that
hypothesis ($H^\prime$) holds is still open.  The main difficulty is that, except for particular cases
of honest times and Jacod's equivalence hypothesis, no explicit formula for the $\gg$-semimartingale
decomposition of an $\ff$-local martingale is available

%%%%%%%%%%%%%%%%%%%%%%%%%%%%%%%%%%%%%%%%%%%%%%%%%%%%%%%%%%%%%%%%%%%%%%
\subsection{Case of a Positive Az\'ema Supermartingale} \lab{sec4.2}
%%%%%%%%%%%%%%%%%%%%%%%%%%%%%%%%%%%%%%%%%%%%%%%%%%%%%%%%%%%%%%%%%%%%%%

\bp \lab{pro1n}
If $Z>0$, then any process $Y^h \in \Classmp $ admits the following representation
\be \lab{xeq3n}
dY^h_t= \left(h_t-\frac{X^h_{t}}{Z_{t}}  \right)
d\mt_t + \frac{1-A_{t-}}{Z_{t-}}\, d\bmz^h_{t} -
(1-A_{t-})\frac{X^h_t}{ Z_t Z_{t-}}\, d\bmz_t \ee
where the first and last terms on the right-hand side are to be understood as Lebesgue--Stieltjes integrals with respect to the finite variation processes $M^{(\tau)}$ and $\overline \mu$.
Equivalently,
upon setting $V^h := X^h Z^{-1}$,
\begin{align} \lab{xeq3na}
dY^h_t &= (h_t-V^h_{t-}) \, d\mt_t + (1-A_{t-}) Z^{-1}_{t-} \,
d\bmz^h_t - (1-A_{t-})V^h_{t-} Z^{-1}_{t-}  \, d\bmz_t \\ & -\Delta
V^h_t \Delta A_t - (1-A_{t-}) Z^{-1}_{t-} \Delta V^h_t \,\Delta
Z_t .  \nonumber
\end{align}
\ep

\proof
We start by noting that for any c\`adl\`ag, $\ff$-adapted process
$V$, the jump process $\Delta V  := V  - V_{ -}$ is a {\it thin}
process, meaning that there exists a sequence $(S_n)_{n\in\mathbb N}$ of
$\ff$-stopping times such that $\{\Delta V \neq 0\}\subset
\bigcup_n [\![S_n]\!]$ (see Definition 7.39 in He et al.
\cite{HWY}). As a consequence, if $B$ is a positive process of
finite variation, then for any c\`adl\`ag and $\ff$-adapted
process $V$ we obtain, for all $t \in \rr_+$,
 \be \lab{exc1}
\int_0^t \Delta V_s \frac{dB_{s}}{B_{s-}} = \sum_{0< s \leq t } \Delta V_s \frac{\Delta B_{s}}{B_{s-}} .
\ee
We denote $X=X^h$ and $V=V^h$ and we proceed as in the proof of Theorem \ref{pox2}. We have that (see \eqref{popo})
\bde
dY^h_t =(h_t- V_{t-}) \, d\myh_t +(1-\myh_{t-})\,
dV_t - \Delta V_t \Delta \myh_t
\ede
and
\begin{align*}
dV_t& = Z^{-1}_{t-}\, dX_t +
X_{t-} \left( - Z^{-2}_{t-}\, dZ_t +Z^{-2}_{t-}\,
\Delta Z_t + \Delta Z^{-1}_t \right)
 + \Delta X_t \Delta Z^{-1}_t .
\end{align*}
Since $dX_t = d\mx_t - h_t \, d\ap_t$ and $dZ_t = d\mu_t- d\ap_t$, we obtain
\bde
Z^{-1}_{t-}\, dX_t - X_{t-}Z^{-2}_{t-}\, dZ_t
 = Z^{-1}_{t-}\, d\mx_t - ( h_t - V_{t-})
 Z^{-1}_{t-} \, d\ap_t - X_{t-} Z^{-2}_{t-}\, d\mu_t,
\ede
which leads to
\be \lab{eq11n}
dY^h_t= (h_t- V_{t-})\, d\mt_t + (1-\myh_{t-}) Z^{-1}_{t-}\, d\mx_t -(1-\myh_{t-}) V_{t-} Z^{-1}_{t-} \, d\mu_t+ dL_t
\ee
where
\begin{align*}
dL_t &:=  (1-\myh_{t-}) \left( V_{t-} Z^{-1}_{t-} \Delta Z_t
 + X_{t-} \Delta Z^{-1}_t  + \Delta X_t \Delta Z^{-1}_t \right)
- \Delta V_t \Delta \myh_t \\ &=-
(1-\myh_{t-}) \Delta V_{t} Z_{t-}^{-1}\Delta Z_t
- \Delta V_t \Delta \myh_t \\
&= -(1-\myh_{t-}) \Delta V_{t} Z_{t-}^{-1} \, dZ_t
 - \Delta V_t \,  d\myh_t = -\Delta V_t \,  d\mt_t
-(1-\myh_{t-}) \Delta V_{t} Z_{t-}^{-1}\, d\mu_t
\end{align*}
where the penultimate equality follows immediately from  \eqref{exc1} and in the
last one we used \eqref{mtau}, the equality $Z = \mu- \ap$ and the definition of $\myh$. The desired integral representation now easily follows from \eqref{eq11n}.  Equality \eqref{xeq3na} is an easy consequence of \eqref{xeq3n}.
\endproof

The integrand in the second term in right-hand side of \eqref{xeq3n} is $\gg$-predictable and
locally bounded so that this term is indeed a $\gg$-local martingale. In contrast, it is not clear a priori whether the first
or the third terms in the right-hand side of \eqref{xeq3n} correspond to $\gg$-local martingales, since the integrands
are not necessarily $\gg$-predictable. The following counterexample shows that the first and the last terms in the right-hand side in \eqref{xeq3n} may indeed fail to be $\gg$-local martingales (although their difference is always a $\gg$-local martingale).

\bex \lab{ex4.1}
We place ourselves within the setup of Proposition \ref{proc1} with $\psi(x)=x$ and $h(x)=x$. Therefore, $h_t = h(N_{t-})
= N_{t-} $ and $\tau = \inf \, \{ t \in \rr_+ \! :\, N_t \geq \Theta \}$ where $\Theta $ has {a} unit exponential
distribution (see \eqref{bbg1}). It is easy to check that the process $h$ belongs to the class $\Classhp$.
Furthermore,  $Z = e^{-N}$ and $X^h = \Phi (N)$ where (see \eqref{nhy})
\bde
 \Phi(x)= \lambda \gamma e^{-x} \, \E \Big( \int_0^\infty (N_s+x) e^{-N_s} \, ds \Big)=e^{-x}(a+bx)
\ede
where we denote $\gamma := 1- e^{-1}$ and where, by elementary computations,
\bde
a := \lambda \gamma \, \E \Big( \int_0^\infty N_s e^{-N_s} \, ds \Big) = \gamma^{-1}-1 , \quad
b := \lambda \gamma \, \E \Big( \int_0^\infty e^{-N_s} \, ds \Big) = 1.
\ede
Consequently, we have that $X^h Z^{-1} =e^{N} \Phi (N) = \wt h (N) = a + N = \gamma^{-1}-1 + N $.
Since the process $\ap$ is absolutely continuous with respect to Lebesgue measure
(see \eqref{mmiy}) and $\wt h(N_s)=\wt h(N_{s-}),\, d\P \otimes ds-$a.e., for the first integral in the right-hand side in \eqref{xeq3n}, we obtain
\begin{align*}
I^1_t & := \int_0^t (h_s - X^h_s Z_s^{-1} ) \, d\mt_s \stackrel{\rm mart}  = \int_0^t  \big( \wt h(N_{s-})  -\wt h(N_{s})  \big) \,d\mt_s \\ &= \int_0^t ( \wt h(N_{s-})  -\wt h(N_{s})  )  \, dA _s - \int_0^{t \wedge \tau } ( \wt h(N_{s-})  -\wt h(N_{s} )  )
 Z_{s-}^{-1} \,dA _s^p  \\ &=  \int_0^t( \wt h(N_{s-})  -\wt h(N_{s})   )  \, dA_s = \int_0^t ( N_{s-}  - N_{s} )  \, dA_s = - A_t
 \end{align*}
where for any two processes $X$ and $Y$, we write $X \stackrel{\rm mart} = Y$ if $X-Y$ is a $\gg$-local martingale and where we used
the property $\{ \Delta A >0 \} \subset \{ \Delta N >0\}$  (see \eqref{ut66}). We conclude that the first term in the right-hand side in \eqref{xeq3n} is not a $\gg$-local martingale.

Let us now examine the last integral in \eqref{xeq3n}. From (\ref{mmiy}), we obtain
\bde
 \mz_t  = 1 + \int_0^t   \big( e^{-\psi (N_{s-}+1)}-e^{-\psi (N_{s-})} \big) \, dM_s = 1 - \int_0^t
\gamma  e^{- N_{s-} }  \, dM_s .
\ede
Due to the immersion property, the process  $M_t = N_t - \lambda t$ is a $\gg$-local martingale and thus
\begin{align*}
&I^3_t  := \int_0^t (1-A_{s-})\frac{X^h_s}{ Z_s Z_{s-}}\, d\bmz_s = \int_0^t (1-A_{s-})\frac{X^h_s}{ Z_s Z_{s-}}\, d\mz_s
=  \int_0^t (1-A_{s-}) Z_{s-}^{-1} \wt h (N_s) \, d\mz_s \\
&=  - \int_0^t (1-A_{s-}) Z_{s-}^{-1} \wt h (N_{s})  \gamma e^{-
N_{s-} } \, dM_s
 \stackrel{\rm mart} = \gamma \int_0^t (1-A_{s-}) (\wt h(N_{s-})-\wt h (N_{s})) \, dM_s
\\ &= \gamma \int_0^t (1-A_{s-}) ( N_{s-} - N_{s}) \, dM_s
= - \gamma \int_0^t (1-A_{s-})\, dM_s - \gamma \int_0^t (1- A_{s-}) \lambda \, ds \stackrel{\rm mart} = - \gamma \lambda (t \wedge \tau )
\end{align*}
so that $I^3$ obviously fails to be a $\gg$-local martingale. However, we also have that
\begin{align*}
I^1_t - I^3_t & \stackrel{\rm mart} = - \big( A_t -  \gamma \lambda (t \wedge \tau ) \big) = - \mt_t
\end{align*}
where the last equality follows from \eqref{mrmiy} applied to $\psi (x)=x$. We conclude that $I^1-I^3$ is indeed a $\gg$-local
martingale, as was expected.

Let us also observe that  \eqref{xeq3n} and the computations above
lead to the following equality
 \begin{align*}
 dY^h_t&= (h_t - \wt h (N_{t-} )  - 1  )\,  d\mt_t + (1-A_{t-})Z^{-1}_{t-} \, d\mx_t   + (1-A_{t-}) \gamma
(\wt h (N_{t-})  + 1)\, dM_t .
\end{align*}
From \eqref{yx1} and \eqref{mmiy}, we obtain
\begin{align*}
\mx_t & = X_t^h + \int_0^t h_s \, d \ap_s = X^h_t + \lambda \int_0^t  N_s \big( e^{- N_{s}}-e^{-N_{s}-1} \big) \, ds
= \Phi (N_t) + \lambda \gamma \int_0^t  N_s e^{- N_{s}} \, ds \\
& = \Phi (0) - \gamma \int_0^t N_{s-}e^{- N_{s-}} \, dN_s + \lambda \gamma \int_0^t N_{s-} e^{- N_{s-}} \, ds
= a -  \gamma \int_0^t N_{s-}e^{- N_{s-}}\, dM_s
\end{align*}
where we also used the relationship $\Phi (x+1) =  \Phi(x) - \gamma x e^{-x}$ (see \eqref{inhy}).
It is clear that $\wt h (x+1) = \wt h(x)+1$ and thus we conclude that
\begin{align*}
 dY^h_t&= (h_t - \wt h (N_{t-} )  - 1  )\,  d\mt_t - (1-A_{t-}) \gamma N_{t-}\, dM_s
     + (1-A_{t-}) \gamma (\wt h (N_{t-})  + 1)\, dM_t \\
&= (h_t - \wt h (N_{t-} +1 )  )\,  d\mt_t  + (1-A_{t-}) \gamma ( - N_{t-} + a + N_{t-} + 1)\, dM_t\\
&= (h_t - \wt h (N_{t-} +1 )  )\,  d\mt_t  + (1-A_{t-})
    (\wt h (N_{t-} +1 ) - \wt h (N_{t-} )) \, dM_t
\end{align*}
since $\gamma (a+1) = 1 = \wt h (N_{t-} +1 ) - \wt h (N_{t-})$.
\eex

\bcor \lab{cor:zpredi}
If $Z>0$ and $\ao=\ap$, then for every process $Y^h \in \Classm $ we
have
\be \lab{xeq3}
dY^h_t= \left(h_t-\frac{X^h_{t}}{Z_{t}}
\right) d\mt_t + \frac{1-A_{t-}}{Z_{t-}}\, d\bmz^h_t
\ee
and the two terms in the right-hand side of \eqref{xeq3} correspond to $\gg$-martingales.
\ecor

\proof
By Lemma \ref{lemma1.1} the postulated equality $\ao=\ap$ implies that $\mu =m=1$ and the representation $Y_t = A_t h_\tau +  (1-A_t) V_t$ holds
and thus the asserted equality \eqref{xeq3} can be deduced from \eqref{xeq3n}.
The $\gg$-martingale property of the first term in
the right-hand side of \eqref{xeq3n} is a consequence of Lemma \ref{l:a3} since indeed $h-Z^{-1}X^h$ is bounded as $|X^h|\leq |\!|h|\!|_\infty Z$.
 The $\gg$-martingale property of the second term in
the right-hand side of \eqref{xeq3n} is clear. \endproof

\bcor \lab{cor1x}
(i) Let $Z>0$, $h\in\Classhp $ and for all $t\in \rr_+$
\be \lab{ass2}
\int_0^t \Delta \mu_s \, d\mt_{s} =  \int_0^t \Delta \mx_s \, d\mt_{s} =  0 .
\ee
Then the process $Y^h \in \Classmp $ admits the following representation
\be  \lab{eqc3}
 dY_t^h= \frac {Z_{t-}}
{Z_{t-}-\Delta \ap_t} \left(h_t-\frac{X^h_{t-}}{Z_{t-}}\right)
d\mt_t+\frac{1-A_{t-}}{Z_{t-}}\, d\bmz^h_t -(1-A_{t-})\frac{X^h_t}{ Z_t Z_{t-}}\, d\bmz_t  .
\ee
(ii) Let $Z>0$, $h\in\Classh$, $\ao=\ap$ and \eqref{ass2} holds. Then
the process $Y^h \in \Classm$ admits the following representation
\be  \lab{beqc3} dY_t^h= \frac {Z_{t-}}
{Z_{t-}-\Delta \ap_t} \left(h_t-\frac{X^h_{t-}}{Z_{t-}}\right)
d\mt_t+\frac{1-A_{t-}}{Z_{t-}}\, d\bmz^h_t . \ee
\ecor

\proof
If condition \eqref{ass2} is satisfied, then from \eqref{xeq3n} and Lemma \ref{lem2}, we obtain
\bde
\left(h_t-\frac{X^h_{t}}{Z_{t}}\right)\, d\mt_t =
\frac{Z_{t-}}{Z_t} \left(h_t-\frac{X^h_{t-}}{Z_{t-}}\right)\, d\mt_t
= \frac {Z_{t-}}  {Z_{t-}-\Delta \ap_t} \left(h_t-\frac{X^h_{t-}}{Z_{t-}}\right)\, d\mt_t
\ede
since $\mt $ is a process of finite variation and $Z_t = Z_{t-}+ \Delta Z_t = Z_{t-} + \Delta \mu_t - \Delta \ap_t$.
Therefore, \eqref{eqc3} is an immediate consequence of \eqref{xeq3n}.
Similarly, the second assertion follows from Corollary \ref{cor:zpredi}.
\endproof

Let us now assume, in addition, that there exists a $d$-dimensional $\ff$-martingale $M$, which has the PRP with respect to the filtration $\ff$. The following result, which is an immediate consequence of Proposition \ref{pro1n}, shows that the $(d+1)$-dimensional $\gg$-martingale $(\mt ,\overline{M})$ generates the class $\Classm $.

\bp  \lab{prox1} Let $Z>0$ and $\ao=\ap$. If an
$\ff$-martingale $M$ has the PRP with respect to $\ff$, then any process $Y^h
\in \Classm $ admits the following representation
\be  \lab{xeq3x}
dY_t^h=\left(h_t-\frac{X^h_{t}}{Z_{t}}\right) \, d\mt_t+\frac{1-A_{t-}}{Z_{t-}}\,\phi^h_t \, d\overline M_t
\ee for some
$\ff$-predictable process $\phi^h $. If, in addition, condition \eqref{ass2} holds, then
\be  \lab{xeq3a}
dY_t^h=\frac
{Z_{t-}}  {Z_{t-}-\Delta \ap_t}
\left(h_t-\frac{X^h_{t-}}{Z_{t-}}\right) d\mt_t +
\frac{1-A_{t-}}{Z_{t-}}\,\phi^h_t \, d\overline M_t .
\ee
\ep

\bex \lab{ex4.2}
Assume that the probability space $(\Omega, \G, \ff,\P)$ supports a random variable $\Theta$  with the unit exponential distribution and such that $\Theta $ is independent of the filtration $\ff$ generated by the Poisson process $N$.
We consider a Cox construction example where
 the random time $\tau$ is given by the expression
\be \lab{bbg}
\tau = \inf \, \{ t \in \rr_+ :\,  \Lambda_t \geq \Theta \}
\ee
where $\Lambda $ is a real-valued, $\ff$-predictable, increasing process such that $\Lambda_0=0$ and $\Lambda_{\infty }:= \lim_{t \to \infty } \Lambda_t = \infty $. Then the immersion property holds for the filtrations $\ff $ and $\gg$ and thus $\tau $ is an $\ff$-pseudo-stopping time. Furthermore, the Az\'ema supermartingale of $\tau $ equals $Z_t = \P (\tau > t \,|\, \F_t) =
e^{-\Lambda_t}$ for all $t \in \rr_+$ and thus it is a strictly positive, decreasing, and $\ff$-predictable process.
We claim then any $\gg$-martingale from the class $\Classm $ admits the integral representation
\be  \lab{xeq3ay}
dY_t^h=\frac{Z_{t-}}  {Z_{t-}-\Delta \ap_t}
\left(h_t-\frac{X^h_{t-}}{Z_{t-}}\right) d\mt_t +
\frac{1-A_{t-}}{Z_{t-}}\,\phi^h_t \, dM_t
\ee
where $M_t := N_t - \lambda t$. Under the present assumptions, we have that $Z = 1- \ap$ and thus $\mu = 1$.  Hence, by Proposition \ref{prox1}, it is enough to show that condition \eqref{ass2} is satisfied by the $\ff$-martingale $\mx $ given by
\eqref{yx2} and satisfying \eqref{mmi} (obviously, $\Delta \mu =0$ so the first equality in \eqref{ass2} is met). In view of definition \eqref{mtau} of $\mt $,  it thus suffices to show that the following equalities are satisfied, for all $t\in\rr_+$,
\be \lab{ass2n}
\int_0^t \Delta \mx_s \, d\ap_{s} =  0 = \int_0^t \Delta \mx_s \, dA_{s}.
\ee
Let us remark that the first equality in \eqref{ass2n} holds if the filtration $\ff$ is quasi-left-continuous and the two equalities
in \eqref{ass2n}  are satisfied when $\tau$ avoids all $\ff$-stopping times. It is well known that the Poisson filtration is quasi-left-continuous and thus the martingale $\mx$ can only jump at $\ff$-totally inaccessible stopping times. Since $\ap$ is $\ff$-predictable, the processes $\mx$ and $\ap $ have no common jumps and thus the first equality in \eqref{ass2n} holds. It now remains to show that the second equality in \eqref{ass2n} is valid for all $t \in \rr_+$. We observe that
\bde
\E \left( \int_0^{\infty} \Delta N_s \, dA_{s} \right) = \int_0^{\infty } \P
( \Delta N_{\tau } = 1 \vert \Theta = \theta ) \, d\P (\Theta \leq
\theta ) = 0
\ede
where the last equality holds since, for any fixed $\theta $, the random time $\tau $ is $\ff$-predictable and
the jump times of the Poisson process are $\ff$-totally inaccessible stopping times. Since $\int_0^{\infty} \Delta N_s \, dA_{s}$ is
non-negative and, in view of \eqref{mmi}, we have that $\{ \Delta \mx \ne 0 \} \subset \{ \Delta M >0\} = \{ \Delta N > 0 \}$,
we conclude that the equality $\int_0^t \Delta \mx_s \, dA_{s} = 0$ is satisfied for
all $t \in \rr_+$. Hence the asserted representation  \eqref{xeq3ay} is a consequence of
equation \eqref{xeq3a} in Proposition \ref{prox1}.
\eex

%%%%%%%%%%%%%%%%%%%%%%%%%%%%%%%%%%%%%%%%%%%%%%%%%%%%%%%%%%%%%%%%%%%%%%%%%%
\subsection{Case of an Arbitrary Az\'ema Supermartingale} \lab{sec4.3}
%%%%%%%%%%%%%%%%%%%%%%%%%%%%%%%%%%%%%%%%%%%%%%%%%%%%%%%%%%%%%%%%%%%%%%%%%%

We relax the assumption that $Z$ is positive and we search for an integral representation of a
$\gg$-martingale $Y^h$ from the class $\Classmp $. In the statement of the proposition, we use the convention
that $0 \cdot \pm \infty = 0$.

\bp \lab{pro1}
Any process $Y^h \in \Classmp $ admits the following
integral representation
\be \lab{eq1111}
dY^h_t=  \left(h_t- \frac{X^h_{t-}}{Z_{t-}} \right)  d\myh_t
-  (1-A_t) \left(h_t-\frac{X^h_{t}}{Z_{t}} \right) \frac{d\ap_t }{Z_{t-}}
+\frac{1-\myh_{t}}{Z_{t-}}\, d\bmz^h_t - (1-\myh_{t})\frac{X^h_{t}}{Z_t Z_{t-}} \, d\bmz_t
\ee
or, equivalently, setting $V^h = X^h Z^{-1}$
\begin{align} \label{eqiu4}
 dY^h_t&=  \left(h_t- V^h_{t-}  \right)
d\mt_t + (1-\myh_{t-})Z^{-1}_{t-}\, d\bmz^h_t -
(1-\myh_{t-}) V^h_{t-} Z^{-1}_{t-} \, d\bmz_t\\
& +   (h_t- V^h_{t-}) \Delta \myh_t Z^{-1}_{t-}\, d\ap_t -
(1-A_{t})\Delta  V^h_t  Z^{-1}_{t-}\, dZ_t
+ \Delta A_t Z^{-1}_{t-}\, (d\bmz^h_t - V^h_{t-} d\bmz_t). \nonumber
\end{align}
\ep

\proof
Let us denote $X= X^h$ and let us set $V := X Z^{-1}$.
Since $Z$ may vanish at time $\tau$, we now introduce the following process
$U_t= (1-A_t)V_t + A_t V_{\tau -}$ so that, in view of \eqref{y}, we obtain $Y_t^h =  (1-A_t) U_t + A_t h_\tau $.
An application of Yoeurp's lemma yields
\begin{equation} \label{xzxz}
dY^h_t =   \left(h_t-\myu_{t-}\right) d\myh_t +(1-\myh_{t})\, d \myu_t
\end{equation}
and, by the integration by parts formula, $dU_t  = (1-\myh_t) \, dV_t$.
The It\^o formula gives, on the set $[\![0,\tau [\![$,
\begin{align*}
dV_t & =  Z^{-1}_{t-}\, dX_t + X_{t-} \left(
- Z^{-2}_{t-}\, dZ_t + Z^{-2}_{t-} \Delta Z_t +
\Delta Z^{-1}_t \right)
 + \Delta X_t \Delta Z^{-1}_t  \\
 &= Z^{-1}_{t-}\, dX_t
- X_{t-}Z^{-2}_{t-}\, dZ_t + \Delta V_t Z^{-1}_{t-} \Delta Z_t
=  Z^{-1}_{t-}\, dX_t
- X_{t-}Z^{-2}_{t-}\, dZ_t + \Delta V_t Z^{-1}_{t-}\, dZ_t \\
&= Z^{-1}_{t-}\, dX_t
-X_{t} (Z_t Z_{t-})^{-1}\, dZ_t
= Z^{-1}_{t-}\, d\mx_t - ( h_t - V_{t}) Z^{-1}_{t-} \, d\ap_t -
 X_{t}(Z_t Z_{t-})^{-1}\,  d\mu_t
\end{align*}
since $dX_t = d\mx_t - h_t \, d\ap_t$ and $dZ_t = d\mu_t- d\ap_t$. By combining the last equality with
\eqref{xzxz}, we obtain \eqref{eq1111}.  Straightforward computations show that equality \eqref{eqiu4} holds as well.
\endproof

\brem
Let us show that if $h\in\Classhp $ and $Z >0$, then representations (\ref{xeq3n}) and  (\ref{eqiu4}) coincide.
On the one hand, by combining the obvious equality $\Delta Y^h_\tau = h_\tau- V_{\tau-}$ with (\ref{xeq3n}), we obtain
\begin{align*} % \lab{eq:jumpsy}
h_\tau- V_{\tau-}&=(h_\tau- V_{\tau})
\big(1- (1-A_{\tau-}) Z^{-1}_{\tau-}\Delta \ap_\tau
\big)+ (1-A_{\tau-}) Z^{-1}_{\tau-} \Delta \mx_\tau -
(1-A_{\tau-}) Z^{-1}_{\tau-} V_\tau \Delta \mu_\tau \nonumber\\&=
h_\tau- V_{\tau} - Z^{-1}_{\tau-} \big( (h_\tau-V_\tau) \Delta \ap_\tau  -
   \Delta \mx_\tau +  V_\tau  \Delta \mu_\tau \big)
\end{align*}
and thus
$$
\Delta V_\tau  + Z^{-1}_{\tau-} \big( (h_\tau-V_\tau) \Delta \ap_\tau  -
\Delta \mx_\tau +  V_\tau  \Delta \mu_\tau \big) =0.
$$
On the other hand, we observe that equation (\ref{eq1111}) can be represented as follows
\begin{align*}
dY^h_t&=(h_t-V_t)\, dM^{(\tau)}_t + (1-\myh_{t-}) Z^{-1}_{t-}\, d\mx_t -
(1-\myh_{t-}) X^h_{t}(Z_t Z_{t-})^{-1} \, d\mu_t
\\& + \Delta V_t \, dA_t+  Z^{-1}_{t-} \Delta A_t \big( (h_t-V_t)\, d\ap_t -
 d\mx_t +  V_t \, d\mu_t \big)
\end{align*}
so that formulae (\ref{xeq3n}) and  (\ref{eqiu4}) coincide.
\erem

\brem After completing this paper, we learnt about the recent work by Choulli et al. \cite{ch} where the
authors introduced the following $\gg$-martingale (see Theorem 2.1 in \cite{ch})
$$
N^{\gg}_t := A_t-\int_0^{t\wedge \tau} \frac{dA^{o,\ff}_s}{\wt Z_{s}}
$$
where $\wt Z_t :=  \P(\tau \geq t \vert \F_t)$. They established several decompositions of a process  $Y^h \in \Classm $
based, in particular, on representations of some $\gg$-martingales as Lebesgue--Stieltjes integrals with respect to
$N^{\gg}$ and with optional integrands (see, in particular, Theorem 2.5 in \cite{ch}).
It is worth pointing out that their method and results require, in principle, a
more detailed information on the conditional distribution of $\tau $ given the reference filtration $\ff$.
 Specifically, the method used in \cite{ch} hinges on the knowledge of
the dual $\ff$-optional projection of $\tau $, whereas our approach relies primarily on the knowledge of the Az\'ema supermartingale $Z$ and its Doob-Meyer decomposition $Z = \mu - \ap $ (hence on the knowledge of the dual $\ff$-predictable projection of $\tau $)
and on the $\gg$-martingale  (see equation \eqref{mtau})
$$
\mt_t := A_t-\int_0^{t\wedge \tau} \frac{dA^{p,\ff}_s}{Z_{s}}.
$$
\erem

% \vskip 5 pt

%%%%%%%%%%%%%%%%%%%%%%%%%%%%%%%%%%%%%%%%%%%%%%%%%%%%%%%%%%%%%%%%%%%%%%%%%%%%%%%%%%%%%%%%%%

\noindent {\bf Acknowledgements.}
The research of Anna Aksamit was supported by the {Chaire  Markets in Transition} (F\'ed\'eration Bancaire Fran\c caise) and by the European Research Council under the European Union's Seventh Framework Programme (FP7/2007-2013)/ERC grant agreement no. 335421.
The research of Monique Jeanblanc was supported by the {Chaire  Markets in Transition  (F\'ed\'eration Bancaire Fran\c caise) and ILB, Labex ANR 11-LABX-0019}. The research of Marek Rutkowski was supported by the DVC Research Bridging Support Grant
{\it Nonlinear Arbitrage Pricing of Multi-Agent Financial Games.} We  warmly thank an anonymous referee for invaluable comments and suggestions, which enabled us to improve our work.

\appendix

\section{Auxiliary Lemmas}

In the proof of Theorem \ref{pox2}, we used the following lemma, which extends a result quoted in Jeulin \cite{Jeu}
(see Remark 4.5 therein). Note that equation \eqref{mtau} can be obtained as a special case of \eqref{mmcc} by setting $\xi = \kappa = 1$.

\bl \label{compensator}
Let the process $B$ be given by the formula $B =\xi A$ where $\xi$ is an integrable and
$\G_\tau$-measurable random variable. Then the process $\wt{M}$ given by the equality
\be \lab{mmcc}
\wt{M}_t = B_t-\int_0^t\frac{1-A_{s-}}{Z_{s-}}\, d\Bpff_s
\ee
is a purely discontinuous $\gg$-martingale stopped at $\tau$. Moreover, the dual $\ff$-predictable projection
of $B$ satisfies $\Bpff_t=\int_0^t \kappa_s \, d\ap _s$ where $\kappa$ is an $\ff$-predictable process such that
the equality $\kappa_\tau=\E (\xi \vert \F_{\tau-})$ holds.
\el

\proof
 Let the process $B$ be given by $B =\xi A$, where the integrable random variable $\xi$ is $\G_\tau$-measurable, and let $\Bpff$ be its dual $\ff$-predictable projection. Recall that $[\![0,\tau ]\!] \subset \{Z_->0\}$ and thus
$$
\wt{M}_t % = B_t-\int_0^t\frac{1-A_{s-}}{Z_{s-}}\, d\Bpff_s
=B_t-\int_0^t\frac{1-A_{s-}}{Z_{s-}}\, \ind_{\{Z_{s-}>0\}} \, d\Bpff_s = B_t-\int_0^t\frac{1-A_{s-}}{Z_{s-}}\, d\wtBpff_s
$$
where we define
$$
\wtBpff_t := \int_0^t \ind_{\{Z_{s-}>0\}} \, d\Bpff_s .
$$
On the one hand, we have, for any $u\geq t$,
\begin{align*}
\E(B_u-B_t\vert\G_t)&=\E(\xi\ind_{\{u\geq \tau>t\}}\vert \G_t)
=\ind_{\{\tau>t\}}Z_t^{-1}\,\E(\xi\ind_{\{u\geq \tau>t\}}\vert\F_t)\\
&= \ind_{\{\tau>t\}}Z_t^{-1}\,\E\left( \int_t^u \xi \ind_{\{Z_{s-}>0\}} \, dA_s \,\Big|\, \F_t\right)
= \ind_{\{\tau>t\}}Z_t^{-1}\,\E\left( \int_t^u \ind_{\{Z_{s-}>0\}} \, dB_s \,\Big|\, \F_t\right)
\\ & = \ind_{\{\tau>t\}}Z_t^{-1}\,\E\left( \int_t^u \ind_{\{Z_{s-}>0\}} \, d\Bpff_s \,\Big|\, \F_t\right) =
\ind_{\{ \tau>t\}}Z_t^{-1}\, \E(\wtBpff_u-\wtBpff_t\vert \F_t).
\end{align*}
On the other hand, we obtain
\begin{align*}
\E \left( \int_t^u\frac{1-A_{s-}}{Z_{s-}}\, d\wtBpff_s \,\Big| \,
\G_t \right)
&=\E \left( \ind_{\{\tau>t\}}\int_t^{u\land\tau}\frac{1}{Z_{s-}}\, d\wtBpff_s \, \Big|\,\G_t \right) \\
& = \ind_{\{\tau>t\}}  Z_t^{-1} \,
\E \left(\ind_{\{\tau>t\}}\int_t^{u\land\tau}\frac{1}{Z_{s-}}\, d\wtBpff_s\,\Big|\,\F_t \right).
\end{align*}
We define the $\ff$-predictable process $\Lambda $ by setting
$\Lambda_t := \int_0^t\frac{1}{Z_{s-}}\, \ind_{\{Z_{s-}>0\}} \, d\Bpff_s  = \int_0^t\frac{1}{Z_{s-}}\, d\wtBpff_s$. Then we obtain, for any $u\geq t$,
\begin{align*}
\E\left(\ind_{\{\tau>t\}}\int_t^{u\land\tau}\frac{1}{Z_{s-}}\, d\wtBpff_s\Vert\F_t\right)
&=\E\left(\ind_{\{\tau>u\}}\int_t^{u}\frac{1}{Z_{s-}}\,\,   d\wtBpff_s+\ind_{\{u\geq\tau>t\}}
\int_t^{\tau}\frac{1}{Z_{s-}}\,\,  d\wtBpff_s\Vert\F_t\right)\\
&=\E\left(Z_u\int_t^{u}\frac{1}{Z_{s-}}\, d\wtBpff_s+\int_t^u\int_t^{v}\frac{1}{Z_{s-}}\, d\wtBpff_s\, d\ap_v\Vert\F_t\right)\\
&=\E\left(Z_u(\Lambda_u-\Lambda_t)+\int_t^u(\Lambda_v-\Lambda_t)\, d\ap_v \Vert\F_t\right)\\
&=\E\left(Z_u(\Lambda_u-\Lambda_t)-\int_t^u(\Lambda_v-\Lambda_t)\, dZ_v \Vert\F_t\right)\\
&=\E \big(\wtBpff_u-\wtBpff_t\vert\F_t \big)
\end{align*}
%We define the $\ff$-predictable process $\Lambda $ by setting $\Lambda_t=\int_0^t\ind_{\{Z_{s-}>0\}}\, d\wtBpff_s$. Then we obtain
%\begin{align*}
%\E\left(\ind_{\{\tau>t\}}\int_t^{u\land\tau}\frac{1}{Z_{s-}}\, d\wtBpff_s\Vert\F_t\right)
%&=\E\left(\ind_{\{\tau>u\}}\int_t^{u}\frac{1}{Z_{s-}}\,\,   d\wtBpff_s+\ind_{\{u\geq\tau>t\}}
%\int_t^{\tau}\frac{1}{Z_{s-}}\,\,  d\wtBpff_s\Vert\F_t\right)\\
%&=\E\left(Z_u\int_t^{u}\mb{\ind_{\{\tau>s\}}\frac{1}{Z_{s-}}}\, d\wtBpff_s+\int_t^u\int_t^{v}\mb{\ind_{\{\tau>s\}}}\frac{1}{Z_{s-}}\, d\wtBpff_s\, d\ap_v\Vert\F_t\right)\\
%&=\E\left(Z_u(\Lambda_u-\Lambda_t)+\int_t^u(\Lambda_v-\Lambda_t)\, d\ap_v \Vert\F_t\right)\\
%&=\E\left(Z_u(\Lambda_u-\Lambda_t)-\int_t^u(\Lambda_v-\Lambda_t)\, dZ_v \Vert\F_t\right)\\
%&=\E \big(\wtBpff_u-\wtBpff_t\vert\F_t \big)
%\end{align*}
where the penultimate equality follows from the fact that for any fixed $t$, the process $(\int_t^u(\Lambda_s-\Lambda_t)\, d\mu_s)_{u\geq t}$ is a true $\ff$-martingale. Indeed, if $\Lambda$ is bounded, then that process is clearly a true $\ff$-martingale and,  since all $\ff$-martingales (in particular, $\mu$) are processes of finite variation, one can use the monotone convergence theorem to prove that $\E(\int_t^u(\Lambda_s-\Lambda_t)\, d\mu_s \vert \F_t)=0$ for all $u\geq t$. The last equality is a consequence of the following computations
\begin{align*}
Z_u&(\Lambda_u-\Lambda_t)-\int_t^u(\Lambda_v-\Lambda_t)\, dZ_v
=Z_u(\Lambda_u-\Lambda_t)+\Lambda_t(Z_u-Z_t)-\int_t^u\Lambda_v \, dZ_v \\
&=Z_u\Lambda_u-\Lambda_tZ_t-\Lambda_uZ_u+\Lambda_tZ_t+\int_t^uZ_{v-}\, d\Lambda_v
=\wtBpff_u-\wtBpff_t.
\end{align*}
The proof of the first statement in Lemma \ref{compensator} is thus completed.
For the second statement, we note that, from the definition of the $\sigma$-field
$\F_{\tau-}$,  the equality $\E(\xi\vert\F_{\tau-})= \kappa_\tau$ holds for some $\ff$-predictable process $\kappa$.
Hence for any bounded, $\ff$-predictable process $X$, we obtain  (note that $X_{\tau }$ is $\F_{\tau-}$-measurable)
\begin{align*}
 \E \left( \int_0^\infty X_s  \, dB^{p,\ff}_s \right) &= \E \left( \int_0^\infty X_s  \, dB_s \right) = \E \left( \int_0^\infty \xi X_s  \, dA_s \right)=
 \E(\xi X_\tau) = \E\big( X_\tau\E(\xi\vert\F_{\tau-})\big) \\& = \E (X_\tau\kappa_\tau)
 = \E \left( \int_0^\infty X_s \kappa_s  \, dA_s \right)= \E \left( \int_0^\infty X_s \kappa_s \, d\ap_s \right)
\end{align*}
where we used the fact that $\ap =A^{p,\ff}$ is the dual $\ff$-predictable projection of $A$. We conclude that $\Bpff_t= \int_0^t
\kappa_s \, d\ap_s$ for all $t \in \rr_+$.
\endproof

In the proof of the next result, we will need the following elementary lemma in which we implicitly assume
that the integrals are well defined.

\bl \lab{lem2}
Let $h\in\Classhp $ and $Z>0$. If $B$ is a process of finite variation such that, for all $t \in \rr_+$,
 \be \lab{assv2}
\int_0^t \Delta \mu_s \, dB_{s} = \int_0^t \Delta \mx_s \, dB_{s} =0
\ee
then  for all $t \in \rr_+$
\bde
\int_0^t\left(h_s-\frac{X^h_{s}}{Z_{s}}  \right) dB_s = \int_0^t
\frac{Z_{s-}}{Z_s} \left(h_s-\frac{X^h_{s-}}{Z_{s-}}\right)\,dB_s.
\ede
\el

\proof
We write $X=X^h$ and $V=XZ^{-1}$. Note that
\bde
\Delta X_t = \Delta \mx_t - h_t  \Delta \ap_t = \Delta \mx_t + h_t  \Delta Z_t - h_t  \Delta \mu_t.
\ede
Using \eqref{assv2}, we obtain
\begin{align*}
(h_t- V_t ) \, dB_t&= (h_t- V_{t-} - \Delta V_t ) \, dB_t \\
& = \left(h_t- V_{t-} -\frac{(X_{t-} + \Delta \mx_t + h_t \Delta Z_t- h_t  \Delta \mu_t ) Z_{t-}- X_{t-}Z_t}{Z_tZ_{t-}} \right) dB_t \\
%& = \frac{h_t Z_t Z_{t-} - X_{t-}Z_t + (X_{t-} + h_t \Delta Z_t) Z_{t-}+ X_{t-}Z_t}{Z_tZ_{t-}}\, dB_t \\
& = \frac{h_t Z_t Z_{t-}  - (X_{t-}+ h_t \Delta Z_t) Z_{t-}}{Z_tZ_{t-}}\, dB_t
 = Z_{t-} Z_t^{-1} (h_t- V_{t-} ) \, dB_t ,
\end{align*}
which gives the desired equality.
\endproof

\bl
\label{l:a3}
Assume that $\mu=1$ so that $Z = 1- \ap $. Then for any $\ff$-optional and bounded process $U$, the process $M^U$
given by the equality
\begin{equation}\label{myx2}
M^U_t =  \int_0^t U_v \, d\mt_v
\end{equation}
is a $\gg$-martingale.
\el

\proof
It is enough to show that $\E ( M^U_u - M^U_t \,|\, \G_t)=0$ for all $u \geq t$. For any $u>t$, we have
\bde
 M^U_u - M^U_t = \int_t^u U_s \, d\mt_s =  U_{\tau } \ind_{\{ t < \tau \leq u \}} - V_{\tau } \ind_{\{ t < \tau \leq u \}} -  V_u \ind_{\{ u < \tau \}}
\ede
where $V$ is the $\ff$-optional process given by $V_u = \int_t^u U_s Z_{s-}^{-1}\ind_{\{Z_{s-}>0\}}\, d\ap_s$ for $s \in [t, \infty )$. Since $Z = 1- \ap $, an application of \eqref{yvv1} and \eqref{yvv3} to $V$ yields
\begin{align*}
&\E ( M^U_u - M^U_t \,|\, \G_t)
= \E\big( U_{\tau } \ind_{\{ t < \tau \leq u \}} \, | \, \G_t \big) - \E\big(  V_{\tau } \ind_{\{ t < \tau \leq u \}} \, | \, \G_t \big) - \E\big(  V_u \ind_{\{ u < \tau  \}} \, | \, \G_t \big)
\\ &= Z^{-1}_t\ind_{\{\tau>t\}} \left [\E\left( \int_t^u U_s \, d\ap_s \, \Big| \, \F_t \right)
- \E\left( \int_t^u V_s \, d\ap_s \, \Big| \, \F_t \right)
-  \E\left( V_u (1 - \ap_u ) \, | \, \F_t \right)\right] = 0
\end{align*}
where the last equality is valid, since the integration by parts formula gives (note that $V_t=0$)
\be \label{vfe}
V_u (1- \ap_u ) + \int_t^u V_s \, d \ap_s= \int_t^u (1- \ap_{s-}) \, dV_s   = \int_t^u  U_s\ind_{\{Z_{s-}>0\}} \, d\ap_s.
\ee
This completes the proof.
\endproof

\bl
\label{l:xz}
(a) For any $\ff$-semimartingale $X$ the process
$\frac{X}{Z}\ind_{\{Z>0\}}\ind_{[\![0,\tau]\!]} $ is a $\gg$-semimartingale and it holds on $[\![0,\tau]\!]$
\begin{align*}
\frac{X_t}{Z_t}\ind_{\{Z_t>0\}}&=\frac{X_0}{Z_0}+\int_0^t\frac{1}{Z_{s-}}\,dX_s
-\int_0^t\frac{X_{s-}}{Z_{s-}^2}\,dZ_s+\sum_{s\leq t}\frac{X_{s-}}{Z_{s-}^2}\, \Delta Z_s
+\sum_{s\leq t}X_s\Delta\bigg(\frac{1}{Z_s}\,\ind_{\{Z_s>0\}}\bigg).
\end{align*}
(b) Let $h$ be a process such that $h_\tau$ is an integrable random variable and let $X^h$ be given by
$$
X^h_t=\E(h_\tau \ind_{\{\tau>t\}} \vert \F_t) = {}^{o,\ff} \big( h_\tau \ind_{[\![ 0, \tau [\![} \big)_t
     = 1- {}^{o,\ff}(h_\tau A) _t.
$$
Then the process $\frac{X^h}{Z}\ind_{\{Z>0\}}\ind_{[\![0,\tau]\!]} $ is a $\gg$-semimartingale and it holds on $[\![0,\tau]\!]$
\begin{align*}
\frac{X^h_t}{Z_t}\ind_{\{Z_t>0\}}&=\frac{X^h_0}{Z_0}+\int_0^t\frac{1}{Z_{s-}}\,dX^h_s
-\int_0^t\frac{X^h_{s-}}{Z_{s-}^2}\, dZ_s
-\sum_{s\leq t}\frac{\Delta Z_s}{Z_{s-}}\, \Delta \bigg(\frac{X^h_s}{Z_s}\,\ind_{\{Z_s>0\}}\bigg).
\end{align*}
\el

\proof
(a) Let us introduce the following $\ff$-stopping times for $n\in \mathbb N$
\begin{align*}
R^n&=\inf \, \{t \in \rr_+ :\, Z_{t}\leq n^{-1}\}, \quad
R=\inf\, \{t \in \rr_+ :\, Z_t=0 \}, \quad \widehat R=R\,\ind_{\{Z_{R-}=0\}}+\infty\ind_{\{Z_{R-}>0\}}.
\end{align*}
Note that $\lim_{n\to\infty}R^n=R$ and $\widehat R$ is an $\ff$-predictable time with the announcing sequence
$(\widehat R^n)_{n\in\mathbb N}$ given by
$$
\widehat R^n=R^n\ \ind_{\{Z_{R^n}>0\}}+\infty\ind_{\{Z_{R^n}=0\}}.
$$
Observe that $\{Z>0\}=[\![0,R[\![$ since $Z$ is a non-negative, right-continuous supermartingale.
Fix $n$.
%By \cite{} \mr{dodac shiqiego i barlowa}
Denote by $Z^-$ the process $Z$ stopped at $R^n-$, i.e., $Z^-_t=Z_t\ind_{\{t<R^n\}}+Z_{R^n-}\ind_{\{t\geq R^n\}}$. Then on $[\![0,R^n]\!]$ we have
\begin{align*}
&\frac{1}{Z_t}\ind_{\{t<R\}}=\frac{1}{Z^-_t}-\ind_{\{t=R^n=R\}}\frac{1}{Z_{R^n-}}+\ind_{\{t=R^n<R\}}\Delta\frac{1}{Z_{R^n}}\\
&=\frac{1}{Z_0}-\int_0^t\frac{1}{(Z^-_{s-})^2}\,dZ^-_s+\sum_{s\leq t}\left (\Delta\frac{1}{Z^-_s}+\frac{1}{(Z^-_{s-})^2}\Delta Z^-_s\right)-\ind_{\{t=R^n=R\}}\frac{1}{Z_{R^n-}}+\ind_{\{t=R^n<R\}}\Delta\frac{1}{Z_{R^n}}\\
&=\frac{1}{Z_0}-\int_0^t\frac{1}{Z_{s-}^2}\,dZ_s+\ind_{\{t=R^n\}}\frac{1}{Z_{R^n-}^2}\Delta Z_{R^n}+\sum_{s\leq t}\left (\Delta\frac{1}{Z^-_s}+\frac{1}{(Z^-_{s-})^2}\Delta Z^-_s\right)\\
&-\ind_{\{t=R^n=R\}}\frac{1}{Z_{R^n-}}+\ind_{\{t=R^n<R\}}\Delta\frac{1}{Z_{R^n}}.
\end{align*}
Note that the last four terms sum up to the common sum and we conclude that on $[\![0,R^n]\!]$ it holds
\begin{align}
\label{eq:z}
&\frac{1}{Z_t}\ind_{\{t<R\}}
=\frac{1}{Z_0}-\int_0^t\frac{1}{Z_{s-}^2}\, dZ_s+\sum_{s\leq t}\left (\Delta\bigg(\frac{1}{Z_s}\ind_{\{s<R\}}\bigg)+\frac{1}{Z_{s-}^2}\Delta Z_s\right).
\end{align}
Therefore, by the integration by parts, on $[\![0, R^n]\!]$ we have that
\begin{align*}
\frac{X_t}{Z_t}\ind_{\{t<R\}}
&=\frac{X_0}{Z_0}+\int_0^t\frac{1}{Z_{s-}}\, dX_s+\int_0^tX_{s-}\, d\left(\frac{1}{Z_s}\ind_{\{s<R\}}\right)+\sum_{s\leq t}\Delta X_s\Delta\bigg(\frac{1}{Z_s}\ind_{\{s<R\}}\bigg)\\
&=\frac{X_0}{Z_0}+\int_0^t\frac{1}{Z_{s-}}\, dX_s
-\int_0^t\frac{X_{s-}}{Z_{s-}^2}\, dZ_s+\sum_{s\leq t}X_{s-}\left (\Delta\bigg(\frac{1}{Z_s}\ind_{\{s<R\}}\bigg)+\frac{1}{Z_{s-}^2}\Delta Z_s\right)\\
&+\sum_{s\leq t}\Delta X_s\Delta\bigg(\frac{1}{Z_s}\ind_{\{s<R\}}\bigg),
\end{align*}
which leads to
\begin{align}
\label{eq:xz}
\frac{X_t}{Z_t}\ind_{\{t<R\}}&=\frac{X_0}{Z_0}+\int_0^t\frac{1}{Z_{s-}}\,dX_s
-\int_0^t\frac{X_{s-}}{Z_{s-}^2}\,dZ_s+\sum_{s\leq t}\frac{X_{s-}}{Z_{s-}^2}\, \Delta Z_s
+\sum_{s\leq t}X_s\Delta\bigg(\frac{1}{Z_s}\ind_{\{s<R\}}\bigg).
\end{align}
Since \eqref{eq:z} and \eqref{eq:xz} hold on $[\![0, R^n]\!]$, they are also valid on $\bigcup_n [\![0, R^n]\!]$.
Note that on $\{\widehat R<\infty\}$
\begin{align*}
0=\Delta Z_{\widehat R}=\Delta m_{\widehat R}-\Delta \ao_{\widehat R}.
\end{align*}
Since $\widehat R$ is an $\ff$-predictable stopping time and $m$ is an $\ff$-martingale, we have $\E(\Delta m_{\widehat R}=0)=0$.
Therefore, by taking expectations and using the definition of the dual optional projection, we obtain
\begin{align*}
0=-\E(\Delta \ao_{\widehat R})=\P(\tau=\widehat R).
\end{align*}
It then follows that $[\![0, \tau]\!]\subset\bigcup_n [\![0, R^n]\!]$ as $\P(\tau\leq R)=1$.
We conclude that identities \eqref{eq:z} and \eqref{eq:xz} hold on $[\![0, \tau]\!]$ and thus the proof of part (a) is complete. \\
(b) Since $X^h$ is an $\ff$-semimartingale, %\mr{optional projection of finite variation process},
 it suffices to apply part (a) and compute the sum terms in this special case. Note that $X^h_R=0$ as $\tau\leq R$ and
$X^h_R=\E(h_\tau \ind_{\{\tau>R\}} \vert \F_R)=0.$ Then we obtain
\begin{align*}
\frac{X^h_{s-}}{Z_{s-}^2}\Delta Z_s
&+X^h_s\Delta\left(\frac{1}{Z_s}\ind_{\{Z_s>0\}}\right)
=\ind_{\{s=R\}}\left( -\frac{X^h_{s-}}{Z_{s-}}-\frac{X^h_s}{Z_{s-}}\right)+\ind_{\{s<R\}}\left( \frac{X^h_{s-}}{Z_{s-}^2}\Delta Z_s+X^h_s\Delta\bigg(\frac{1}{Z_s}\bigg)\right)\\
&=-\ind_{\{s=R\}}\frac{X^h_{s-}}{Z_{s-}}-\ind_{\{s<R\}} \frac{\Delta Z_s}{Z_{s-}}\,\Delta \bigg(\frac{X^h_s}{Z_s}\bigg)
=-\frac{\Delta Z_s}{Z_{s-}}\Delta \bigg(\frac{X^h_s}{Z_s}\ind_{\{Z_s>0\}}\bigg)
\end{align*}
since $X^h_R=0$ and $-\frac{\Delta Z_R}{Z_{R-}}=1$.
\endproof

%%%%%%%%%%%%%%%%%%%%%%%%%%%%%%%%%%%%%%%%%%%%%%%%%%%%%%%%%%%%%%%%%%%%%%%%%%%%%%%%%%%%%%%

\end{document}